\documentclass[final,onefignum,onetabnum]{siamart171218}

\usepackage{lipsum}
\usepackage{amsfonts}
\usepackage{graphicx}
\usepackage{epstopdf}
\usepackage{algorithmic}
\usepackage{bm}
\ifpdf
  \DeclareGraphicsExtensions{.eps,.pdf,.png,.jpg}
\else
  \DeclareGraphicsExtensions{.eps}
\fi

\newsiamremark{remark}{Remark}
\newsiamremark{hypothesis}{Hypothesis}
\crefname{hypothesis}{Hypothesis}{Hypotheses}
\newsiamthm{claim}{Claim}

\headers{Stability Analysis of QBMM for Kinetic Equations}{Q. Huang, S.Q. Li, W.A. Yong}

\title{Stability Analysis of Quadrature-based Moment Methods for Kinetic Equations\thanks{Submitted to the editors DATE.
\funding{This work was funded by China Postdoctoral Science Foundation (under contract no.~043201001) and National Natural Science Foundation of China (no.~51725601 and 11471185).}}}

\author{Qian Huang\thanks{Zhou Pei-Yuan Center for Applied Mathematics, Tsinghua University, Beijing 100084, China
  (\email{hqqh91@qq.com}, \url{https://www.researchgate.net/profile/Qian_Huang34}).}
\and Shuiqing Li\thanks{Department of Energy and Power Engineering, Tsinghua University, Beijing 100084, China
  (\email{lishuiqing@tsinghua.edu.cn}, \url{http://www.thu-lishuiqing.org}).}
\and Wen-An Yong\thanks{Department of Mathematical Sciences, Tsinghua University, Beijing 100084, China (\email{wayong@tsinghua.edu.cn}, \url{https://www.researchgate.net/profile/Wen-An_Yong}).}}

\usepackage{amsopn}
\DeclareMathOperator{\diag}{diag}

\ifpdf
\hypersetup{
  pdftitle={Stability Analysis of Quadrature-based Moment Methods for Kinetic Equations}
  pdfauthor={Q. Huang, S.Q. Li, W.A. Yong}
}
\fi

\begin{document}

\maketitle

\begin{abstract}
    In this paper, we present a systematic stability analysis of the quadrature-based moment method (QBMM) for the one-dimensional Boltzmann equation with BGK or Shakhov models. As reported in recent literature, the method has revealed its potential for modeling non-equilibrium flows, while a thorough theoretical analysis is largely missing but desirable. We show that the method can yield non-hyperbolic moment systems if the distribution function is approximated by a linear combination of $\delta$-functions.
    On the other hand, if the $\delta$-functions are replaced by their Gaussian approximations with a common variance, we prove that the moment systems are strictly hyperbolic and preserve the dissipation property (or $H$-theorem) of the kinetic equation.
    In the proof we also determine the equilibrium manifold that lies on the boundary of the state space.
    The proofs are quite technical and involve detailed analyses of the characteristic polynomials of the coefficient matrices.
\end{abstract}



\section{Introduction}

Kinetic theories pioneered by L. Boltzmann arise in a variety of fields beyond the classical rarefield gas dynamics, ranging from multiphase flows \cite{HQ2017,MarFox2013}, aerosol dynamics in atmospheric environments \cite{Frd2000,MarFox2013}, and active matter physics \cite{Herty2018}, to galactic dynamics in the universe \cite{Yoshi2013}. In the kinetic framework \cite{Kremer2010}, various physical systems are described with a distribution function $f$ which depends on the spatial and other problem-specific microscopic variables and its time evolution is governed by kinetic equations like the Boltzmann equation. Although the kinetic equations have solid physical ground, they are computationally costly and therefore not directly usable in engineering applications.

Because of the above reason, various simplifications or approximations of the kinetic equations have been proposed, including the BGK model \cite{bgk1954}, discrete velocity models \cite{Gatignol1975,Mieu2000}, and moment closure systems \cite{Grad1949,Levermore1996,MarFox2013}. All these approximations have their advantages and disadvantages. This work is concerned with moment closure systems, in which the governing equations of several moments of the distribution function are derived from the kinetic equation and an additional procedure must be accompanied to close the moment system \cite{MarFox2013}. The resultant moment systems consist usually of first-order partial differential equations (PDEs).

To correctly model the observability of physical processes, the derived system of PDEs should be well-posed (or hyperbolic for first-order systems). For instance, the well-known Grad's closure method yields non-hyperbolic PDEs and produces unphysical results \cite{Grad1949,Muller1998}. Its hyperbolic regularization has attracted much attention \cite{Cai13-1,Cai13-2,Levermore1996,McDonald2013,Struchtrup2003}. A recent work is \cite{Cai2015} where the authors introduced a framework to construct hyperbolic moment closure systems.

Furthermore, the moment closure systems derived from the kinetic equations should preserve the key physical properties of the original kinetic equations. For the Boltzmann equation, one of the key properties is the celebrated $H$-theorem characterizing the dissipation property of the mesoscopic system under consideration \cite{Kremer2010}. In this regard, a paradigm is the widely used BGK model that not only simplifies the collision term in the Boltzmann equation, but also inherits the key conservation and dissipation properties thereof \cite{Kremer2010}. At this point, an immediate question is how to manifest the $H$-theorem in such moment systems.

It turns out that the structural stability condition proposed in \cite{Yong1999} for hyperbolic relaxation systems is a proper counterpart of the $H$-theorem for the kinetic equation. Indeed, this condition has been tacitly respected by many well-developed physical theories \cite{Yong2008}. Recently, it was shown in \cite{Di2017} to be satisfied by the hyperbolic regularization models derived in \cite{Cai13-1,Cai13-2,Cai2015}.
In contrast, the Biot/squirt (BISQ) model for wave propagation in saturated porous media violates this condition and thus allows exponentially exploding asymptotic solutions \cite{LiuJW2016}. On the other hand, this condition also implies that the resultant moment system is compatible with the classical theories \cite{Yong1999}. The implication is important because the lower-order moments are usually associated with the macroscopic parameters of the system \cite{Kremer2010}. Therefore, we believe that the structural stability condition is a proper criterion to evaluate the moment closure systems.

The objective of this paper is to investigate whether or not the quadrature-based moment method (QBMM, \cite{MarFox2013}) yields hyperbolic PDEs which satisfy the structural stability condition above.
In QBMM, the distribution function $f$ is approximated with a linear combination of $N$ ($N \ge 1$) $\delta$-functions with unknown centers or their Gaussian approximations with unknown variance and centers (named QMOM or EQMOM, respectively) \cite{MarFox2013}.
QBMM has become an effective and popular method in simulating the evolution of fine particulate matter, where the distribution function is independent of the particle velocity and the resultant governing equation is termed population balance equation \cite{MarFox2013,Ngu2016,Yuan2011}.
However, the QMOM-derived moment system of the Boltzmann equation leads to unphysical shocks in the numerical solution of Riemann problems \cite{Fox2008}, which is confirmed by our own numerical results (see the Supplementary Material).
Thus it is appealing to find the cause for the irregular behaviors and the aforementioned criteria are expected to be useful in clarifying such issues.

This paper deals only with the spatial one-dimensional (1-D) Boltzmann equation with hypothetical collisions (BGK or Shakhov type), just to figure out a road map for further investigations of general cases. We show that the QMOM-derived moment system is not strongly hyperbolic for any number $N$ of nodes, while the Gaussian EQMOM produces strictly hyperbolic moment systems when the variance is positive. For the latter, we further determine their equilibrium manifolds and verify the structural stability condition. The proofs are quite technical and purely analytic. They involve detailed analyses of characteristic polynomials of the coefficient matrices.

Let us remark that for $N=2$, the hyperbolicity of moment systems has been studied in \cite{Chalons2012} for 1-D QMOM and in \cite{Chalons2017} for 1-D Gaussian-EQMOM. The proofs rely on direct calculations of the eigenvalues of the coefficient matrix of the moment systems \cite{Chalons2012,Chalons2017} and
does not seem generalizable to $N$-node systems. Thus new techniques are needed to handle the general cases. Moreover, the stability of EQMOM has not been analyzed in the existing literature. Given our positive results, EQMOM reveals its potential in solving a wider range of kinetic equations.

The paper is organized as follows. \Cref{sec:main} presents a brief introduction on QBMM (QMOM and EQMOM) and states our main results. \Cref{sec:StabQMOM} is devoted to a proof of non-hyperbolicity of QMOM for $N$-node systems.
In \Cref{sec:StabEQMOM}, we verify the structural stability condition for the EQMOM with $N$ nodes. In particular, the hyperbolicity is demonstrated in \Cref{subsec:hypeqmom}, the equilibrium states are determined in \Cref{subsec:eqmomequistate}, and the dissipation property is shown in \Cref{subsec:coupAseqmom,subsec:coupASseqmom}. Finally, we conclude our paper in \Cref{sec:conclusions}.

\section{Preliminaries}
\label{sec:main}

For simplicity, we only consider a hypothetical 1-D ideal gas with the probability density function $f=f(t,x,\xi)$ of time $t \in \mathbb{R}^{+}$, spatial position $x \in \mathbb{R}$ and velocity $\xi \in \mathbb{R}$.
The temporal evolution of $f$ is governed by the Boltzmann equation \cite{Kremer2010}:
\begin{equation} \label{eq:boltz}
    \frac{\partial f}{\partial t} + \xi \frac{\partial f}{\partial x} = Q(f).
\end{equation}
Here the volumetric force is neglected and the right-hand side $Q(f)$ represents the collisions. As a standard assumption \cite{Kremer2010}, $Q=Q(f)$ has only 1, $\xi$ and $\xi^2$ as locally conserved quantities:
\begin{equation} \label{eq:collinvar}
    \int_{\mathbb{R}} Q(f) \phi(\xi) d \xi = 0, \qquad \phi (\xi) = 1, \ \xi, \ \xi^{2},
\end{equation}
and vanishes at a local equilibrium distribution
\begin{equation} \label{eq:feq}
    f_{eq} = f_{eq} (t,x,\xi) = \frac{\rho}{(2\pi \theta)^{1/2}} \exp \left( -\frac{ \left( \xi - U \right)^{2} }{2 \theta} \right),
\end{equation}
where $\rho$, $U$ and $\theta$ are the density, velocity and temperature of the gas, respectively. They are the classical macroscopic parameters related to $f$ as
\begin{equation} \label{eq:mo2macro}
    \rho = \int_{\mathbb{R}} f d \xi, \quad
    \rho U = \int_{\mathbb{R}} \xi f d \xi, \quad
    \rho \theta = \int_{\mathbb{R}} \left( \xi-U \right)^{2}  f d \xi,
\end{equation}

In this paper, we mainly consider the BGK model \cite{bgk1954}, where
\begin{equation} \label{eq:BGK}
    Q=Q_{BGK}(f) = \nu (f_{eq} - f).
\end{equation}
Here $\nu$ is the collision frequency.
This simple model has been widely used since it preserves several key properties of the kinetic equation, including \cref{eq:collinvar} and the $H$-theorem. Because the BGK model results in the Prantle number $Pr=1$, inconsistent with most realistic cases \cite{Kremer2010}, the Shakhov model was proposed \cite{Shakhov1968}:
\begin{equation} \label{eq:Shakhov}
    Q_{S}(f) = \nu (f_S - f).
\end{equation}
Here an alternative equilibrium distribution $f_S$ is assumed:
\begin{equation} \label{eq:fS}
    f_S = f_{eq} \times \left( 1 + \frac{(1-Pr) q (\xi-U)}{3 \rho \theta^2} \left( \frac{ \left( \xi-U \right)^{2}}{\theta} - 3 \right) \right)
\end{equation}
with $q$ the heat flux defined as $q = \int_{\mathbb{R}} \frac{1}{2} \left( \xi-U \right)^3  f d \xi$.

Denote by $M_j(t,x) = \int_{\mathbb{R}} \xi^{j} f d\xi$ the $j$th velocity-moment of $f$.
From \cref{eq:mo2macro} we see that
\begin{equation} \label{eq:1Dmo}
    M_0 = \rho, \ M_1 = \rho U, \ M_2 = \rho (U^2+\theta), \ M_3 = \rho (U^3+3U\theta)+2q.
\end{equation}
The evolution equation for $M_j$ can be derived from the Boltzmann equation \cref{eq:boltz} with the BGK collision \cref{eq:BGK}:
\begin{equation} \label{eq:1Dmoeqbgk}
    \partial_t M_j + \partial_x M_{j+1} = \nu \left[ \rho \Delta_j(U,\theta) - M_j \right].
\end{equation}
Here $\Delta_j(U,\theta)$ denotes the $j$th moment of the normalized Gaussian distribution
\begin{displaymath}
    \delta_{\theta} (\xi;U) = \frac{1}{\sqrt{2\pi \theta}} \exp \left( - \frac{(\xi-U)^2}{2\theta} \right).
\end{displaymath}
Notice that $\Delta_0(U,\theta)=1$ and $\Delta_1(U,\theta)=U$.

There are infinitely many equations in \cref{eq:1Dmoeqbgk}. The first $N$ equations for moments $M_0,\dots,M_{N-1}$ are not closed, because the $M_{N-1}$-equation contains the term $\partial_x M_{N}$. Hence a closure method is needed.

In the rest of this section, we introduce the QBMM methods, the structural stability condition for hyperbolic relaxation systems, and our main results of this paper.

\subsection{Quadrature-based moment methods}
\label{subsec:QBMM}

In QBMM, the lower-order moments determine the weights and nodes of the quadrature for the integration $\int f(\xi) g(\xi) d \xi$.
Then the unclosed term can be expressed in terms of the lower-order moments and thereby the closure is done \cite{MarFox2013}.

\subsubsection{Quadrature method of moment (QMOM)}
\label{subsubsec:QMOM}

In QMOM, the distribution function $f$ is assumed to be a sum of $N$ Dirac delta functions
\begin{equation} \label{eq:QMOMf}
    f(\xi) = \sum_{i=1}^{N} w_i \delta(\xi-u_{i}).
\end{equation}
In order to determine the weights $w_i$ and nodes $u_i$, the first $2N$ lower-order moments $M_0,\dots,M_{2N-1}$ are employed:
\begin{equation} \label{eq:QMOMmo}
    M_j = \sum_{i=1}^{N} w_i u_i^j \quad \text{for } j=0,\dots,2N-1.
\end{equation}
These non-linear algebraic equations can be solved to obtain $w_i$ and $u_i$ as in \cite{MarFox2013}. Then the next moment $M_{2N}$ can be found as
\begin{equation} \label{eq:QMOMM_2N}
    \bar{M}_{2N} = \sum_{i=1}^{N} w_i u_i^{2N}.
\end{equation}
Namely, $w_i$ and $u_i$ are functions of $M_1,\dots,M_{2N-1}$, and so is $\bar{M}_{2N}$. In this way, we obtain the following system of PDEs:
\begin{equation} \label{eq:QMOMboltz}
    \partial_t M + A(M) \partial_x M = \nu S(M).
\end{equation}
Here $M=(M_0,\dots,M_{2N-1})^{T} \in \mathbb{R}^{2N}$, $S(M) = \rho (\Delta_0 (U,\theta),\dots, \Delta_{2N-1} (U,\theta))^{T} - M$, and
\begin{equation} \label{eq:QMOMA}
    A(M) =
    \begin{bmatrix}
        0 & 1 & & & \\
        & 0 & 1 & & \\
        & & \ddots & \ddots & \\
        & & & 0 & 1 \\
        a_0 & a_1 & \cdots & a_{2N-2} & a_{2N-1}
    \end{bmatrix},
\end{equation}
with $a_j = \frac{\partial \bar{M}_{2N}}{\partial M_j}$ for $0 \le j \le 2N-1$.

\subsubsection{Extended-QMOM (EQMOM)}
\label{subsubsec:EQMOM}

In order to improve QMOM \cite{Chalons2017}, the delta function in \cref{eq:QMOMf} is replaced with its Gaussian approximation
\begin{displaymath}
    \delta_{\sigma^2} (\xi;u) = \frac{1}{\sqrt{2\pi} \sigma} \exp \left( - \frac{(\xi-u)^2}{2\sigma^2} \right),
\end{displaymath}
that is,
\begin{equation} \label{eq:EQMOMf}
    f(\xi) = \sum_{i=1}^{N} w_i \delta_{\sigma^2}(\xi;u_{i}).
\end{equation}
Set $W = (w_1,u_1,\dots,w_N,u_N,\sigma^2)^{T} \in \mathbb{R}^{2N+1}$ and $M=(M_0,\dots,M_{2N})^{T} \in \mathbb{R}^{2N+1}$. They are related with the map $M = \mathcal{M} (W)$:
\begin{equation} \label{eq:EQMOMmo}
    M_j = \sum_{i=1}^{N} w_i \Delta_j(u_i,\sigma^2) \quad \text{for } j=0,\dots,2N,
\end{equation}
defined for $W \in \Omega_W = \Omega_W^{open} \cup \Omega_{W}^{eq}$, where
\begin{subequations} \label{eq:domainW}
\begin{align}
\Omega_{W}^{open} &= \{ W: \ w_i>0; \ \sigma^2>0; \ \forall \ i \ne j, \ u_i \ne u_j \}, \label{eq:domainWopen} \\
\Omega_{W}^{eq} &= \{ W: \ w_i>0; \ \sigma^2 > 0; \ u_1=u_2=\dots=u_N \}.
\end{align}
\end{subequations}

Remark that $\Delta_j(u_i, \sigma^2)$ is exactly the same as that in \cref{eq:1Dmoeqbgk}. It is shown in \Cref{sec:injection} that the map $M = \mathcal{M} (W)$ is one-to-one for $W \in \Omega_{W}^{open}$.
Therefore, $W$ can be uniquely solved from \cref{eq:EQMOMmo} for $M \in \mathcal{M}(\Omega_W^{open})$.
In this way, the next moment $M_{2N+1}$ is a function of the lower-order moments $M \in \mathcal{M}(\Omega_W^{open})$:
\begin{equation} \label{eq:EQMOMM_2Np1}
    \bar{M}_{2N+1} = \sum_{i=1}^{N} w_i \Delta_{2N+1} (u_i, \sigma^2).
\end{equation}
Therefore, the following moment system is derived:
\begin{equation} \label{eq:EQMOMboltz}
    \partial_t M + A(M) \partial_x M = \nu S(M)
\end{equation}
 for $M \in \mathcal{M} (\Omega_{W}^{open})$.
Here $S(M) = \rho \left(\Delta_0 (U,\theta),\dots, \Delta_{2N} (U,\theta) \right)^{T} - M$ and
\begin{equation} \label{eq:EQMOMA}
    A(M) =
    \begin{bmatrix}
        0 & 1 & & & \\
        & 0 & 1 & & \\
        & & \ddots & \ddots & \\
        & & & 0 & 1 \\
        a_0 & a_1 & \cdots & a_{2N-1} & a_{2N}
    \end{bmatrix}
\end{equation}
with $a_j = \frac{\partial \bar{M}_{2N+1}}{\partial M_j}$ for $0 \le j \le 2N$.

For such systems, the moment set $\mathcal{M} (\Omega_{W})$ and its closure $\overline{\mathcal{M} (\Omega_{W})}$ have been extensively studied as a realizability issue in the literature \cite{Chalons2017,MarFox2013,Ngu2016,Pigou2018}.
A further discussion on this issue is beyond the scope of this paper.

\subsection{Structural stability condition}
\label{subsec:stability}

The both QMOM and EQMOM moment systems consist of first-order PDEs derived from the Boltzmann equation. To clarify whether or not these systems inherits the $H$-theorem characterizing the dissipation property of the Boltzmann equation, we recall the structural stability condition proposed in \cite{Yong1999} for systems of $D$-dimensional PDEs:
\begin{equation} \label{eq:genehypde}
    \frac{\partial M}{\partial t} + \sum_{d=1}^{D} A_d (M) \frac{\partial M}{\partial x_d} = S(M).
\end{equation}
Here $M$ is the unknown $n$-vector valued function, $A_d=A_d(M)$ is the $d$th $n \times n$ coefficient matrix, and the source term $S=S(M)$ is a given $n$-vector valued function of $M \in \mathbb{G} \subset \mathbb{R}^n$. As in \cite{Yong1999}, we assume that the equilibrium manifold $\mathcal{E} = \{ M \in \mathbb{G} \ | \ S(M) = 0 \}$ is not empty and denote the Jacobian matrix of $S(M)$ as $S_{M} (M)$. The stability condition reads as
\begin{itemize}
    \item [(i)] There exist an invertible $n \times n$ matrix $P(M)$ and an invertible $r \times r$ ($0<r \le n$) matrix $\hat{T} (M)$ such that
    \begin{displaymath}
        P(M) S_M (M) =
        \begin{bmatrix}
            0 & 0 \\
            0 & \hat{T}(M)
        \end{bmatrix}
        P(M), \quad \forall M \in \mathcal{E};
    \end{displaymath}

    \item [(ii)] There exists a positive definite symmetric matrix $A_0(M)$ such that
    \begin{displaymath}
        A_0 (M) A_d (M) = A_d^T (M) A_0 (M) \quad \text{for any } M \in \mathbb{G} \text{ and } d = 1,\dots,D;
    \end{displaymath}

    \item [(iii)] The spatial derivative parts and the source are coupled as
    \begin{displaymath}
        A_0 (M) S_M (M) + S_M^T (M) A_0 (M) \le - P^T (M)
        \begin{bmatrix}
            0 & 0 \\
            0 & I_r
        \end{bmatrix}
        P(M), \quad \forall M \in \mathcal{E}.
    \end{displaymath}
\end{itemize}
Here $I_r$ is the unit matrix of order $r$.

As shown in \cite{Yong2008}, this set of conditions has been tacitly respected by many well-developed physical theories. Condition (i) is classical for initial value problems of the system of ordinary differential equations (ODE, spatially homogeneous systems), while (ii) means the symmetrizable hyperbolicity of the PDE system. Condition (iii) characterizes a kind of coupling between the ODE and PDE parts. Recently, this structural stability condition is shown in \cite{Di2017} to be proper for certain moment closure systems.
On the other hand, this set of conditions implies the existence and stability of the zero relaxation limit of the corresponding initial value problems \cite{Yong1999}. Thanks to these, we believe that the structural stability condition is essential for a reasonable moment closure system.

\subsection{Main results}
\label{subsec:mainresult}

For the moment systems derived above, we will establish the following facts as the main result of this paper,

\begin{theorem}[Non-hyperbolicity of QMOM] \label{thm:stabqmom}
    The QMOM-derived moment system \cref{eq:QMOMboltz} is not strongly hyperbolic.
\end{theorem}

\begin{theorem}[Stability of EQMOM] \label{thm:stabeqmom}
    The EQMOM-derived moment system \cref{eq:EQMOMboltz} satisfies the structural stability condition for $M \in \mathcal{M}(\Omega_W)$.
\end{theorem}

A proof of \cref{thm:stabqmom} will be presented in the next section. In  \Cref{sec:StabEQMOM}, \cref{thm:stabeqmom} is divided as \cref{thm:hypeqmom} (hyperbolicity of EQMOM), \cref{thm:equisolu} (equilibrium state), \cref{thm:eqmombgkstabiii} (BGK model) and \cref{thm:eqmomshakstab} (Shakhov model), which will be proved in \Cref{subsec:hypeqmom,subsec:eqmomequistate,subsec:coupAseqmom,subsec:coupASseqmom}, respectively.

\section{Non-hyperbolicity of QMOM}
\label{sec:StabQMOM}

This section is devoted to a proof of \cref{thm:stabqmom} for the QMOM-derived moment system \cref{eq:QMOMboltz} with $N \ge 2$. We should mention that this theorem has been proved in \cite{Chalons2012} but only for $N=2$. For our purpose, we need to consider the  $2N \times 2N$ coefficient matrix $A=A(M)$ in \cref{eq:QMOMA}.

\begin{proof}[Proof of \cref{thm:stabqmom}]
    Let $\lambda$ be an eigenvalue of $A$ and $\mathbf{v}=(v_1,\dots,v_{2N})^{T}$ the corresponding right eigenvector.
    A direct calculation indicates that
    \begin{subequations}
    \begin{align}
        v_{k} = \lambda v_{k-1} &= \lambda^{k-1} v_1 \quad \text{for } k=2,\dots,2N, \label{eq:Avseq} \\
        \sum_{k=1}^{2N} a_{k-1} v_k &= \lambda v_{2N} = \lambda^{2N} v_1. \label{eq:Avlast}
    \end{align}
    \end{subequations}
    Then we have $\mathbf{v} = v_1 (1,\lambda,\dots,\lambda^{2N-1})^{T}$ and thereby $v_1 \ne 0$. This shows that the geometric multiplicity of each eigenvalue is 1.

    On the other hand, we see from \cref{eq:Avlast} that the characteristic polynomial of $A$ is
    \begin{equation} \label{eq:Acharpqmom}
        c(\lambda) = \lambda^{2N} - a_{2N-1} \lambda^{2N-1} - \dots - a_1 \lambda - a_0.
    \end{equation}
    Note that $(a_0,a_1,\dots,a_{2N-1}) = \left( \frac{ \partial \bar{M}_{2N}}{\partial M_0}, \frac{\partial \bar{M}_{2N}}{\partial M_1}, \dots, \frac{\partial \bar{M}_{2N}}{\partial M_{2N-1}} \right) = \frac{\partial \bar{M}_{2N}}{\partial M}$ with $\bar{M}_{2N}$ defined in \cref{eq:QMOMM_2N} and $M=(M_0,\dots,M_{2N-1})^{T}$.
    Writing $W = (w_1,u_1,\dots,w_N,u_N)^{T} \in \mathbb{R}^{2N}$, we have
    \begin{equation} \label{eq:aMWqmom}
        (a_0,a_1,...,a_{2N-1}) \frac{\partial M}{\partial W} = \left( \frac{\partial \bar{M}_{2N}}{\partial M} \right) \left( \frac{\partial M}{\partial W} \right) = \frac{\partial \bar{M}_{2N}}{\partial W}.
    \end{equation}
    In addition, it follows from \cref{eq:QMOMmo} that the $2N \times 2N$ Jacobian matrix $\partial M / \partial W$ is
    \begin{displaymath}
    \frac{\partial M}{\partial W} =
    \begin{bmatrix}
        1 & 0 & \cdots & 1 & 0 \\
        u_1 & w_1 & \cdots & u_N & w_N \\
        \vdots & \vdots & & \vdots & \vdots \\
        u_1^j & jw_1u_1^{j-1} & \cdots & u_N^j & jw_Nu_N^{j-1} \\
        \vdots & \vdots & & \vdots & \vdots \\
        u_1^{2N-1} & (2N-1)w_1u_1^{2N-2} & \cdots & u_N^{2N-1} & (2N-1)w_Nu_N^{2N-2}
    \end{bmatrix}
    \end{displaymath}
    and from \cref{eq:QMOMM_2N} that
    \begin{displaymath}
        \frac{\partial \bar{M}_{2N}}{\partial W} = \left( u_1^{2N}, 2Nw_1u_1^{2N-1},\dots, u_N^{2N}, 2Nw_Nu_N^{2N-1} \right).
    \end{displaymath}
    Substituting the last two relations into \cref{eq:aMWqmom}, we obtain
    \begin{displaymath}
    \begin{aligned}
        u_k^{2N} - a_{2N-1} u_k^{2N-1} - \dots - a_1 u_k - a_0 &= 0, \\
        2Nu_k^{2N-1} - (2N-1)a_{2N-1}u_k^{2N-2} - \dots - a_1 &= 0
    \end{aligned}
    \end{displaymath}
    for $k=1,\dots,N$. These mean that $c(u_k) = 0$ and $ \left. \frac{dc(\lambda)}{d \lambda} \right|_{\lambda = u_k} = 0$ for $k=1,\dots,N$. Since $c=c(\lambda)$ is a monic polynomial of order $2N$, there must be
    \begin{equation} \label{eq:cuk_qmom}
        c(\lambda) = (\lambda-u_1)^2 \cdots (\lambda-u_N)^2.
    \end{equation}

    As a result, the eigenvalues of $A$ are $u_1,u_2,\dots,u_N$ and each of them has the algebraic multiplicity 2 and the geometric multiplicity 1. In view of its Jordan canonical form, the coefficient matrix $A$ is similar to
    \begin{equation} \label{eq:AJor_qmom}
    \begin{bmatrix}
        u_1 & 1 & & & \\
        0 & u_1 & & & \\
         & & \ddots & & \\
         & & & u_N & 1 \\
         & & & 0 & u_N
    \end{bmatrix}.
    \end{equation}
    Hence the moment closure system \cref{eq:QMOMboltz} is not strongly hyperbolic.
\end{proof}

\section{Stability of EQMOM}
\label{sec:StabEQMOM}

We prove \cref{thm:stabeqmom} in this section. In particular,  \Cref{subsec:hypeqmom} is devoted to Condition (ii), while Conditions (i) and (iii) are verified in \Cref{subsec:coupAseqmom,subsec:coupASseqmom} for both the BGK and Shakhov collision models.

\subsection{Preliminaries}
\label{subsec:preli}

Recall that in \Cref{subsec:QBMM}, we use the notation
\begin{displaymath}
    \Delta_{j} = \Delta_{j}(u,\sigma^2) = \int_{\mathbb{R}} \xi^j \delta_{\sigma^2} (\xi;u) d\xi
\end{displaymath}
for the $j$th moment of the Gaussian distribution $\delta_{\sigma^2}=\delta_{\sigma^2}(\xi;u) = \frac{1}{\sqrt{2\pi}\sigma} \exp{\left(-\frac{(\xi-u)^2}{2\sigma^2}\right)}$.
A direct calculation shows $\Delta_0(u,\sigma^2) = 1$ and $\Delta_1(u,\sigma^2) = u$. Moreover, we can show with \cref{prop:basicDel}(a) below that $\Delta_j (u,\sigma^2)$ is a bivariate polynomial of $u$ and $\sigma^2$.

\begin{lemma} \label{prop:basicDel}
    \begin{displaymath}
    \begin{aligned}
        &\text{(a)}& \ \Delta_j(u,\sigma^2) &= u \Delta_{j-1}(u,\sigma^2) + (j-1) \sigma^2 \Delta_{j-2}(u,\sigma^2) \quad \text{for } j \ge 2, \\
        &\text{(b)}& \ \Delta_j(u,\sigma^2) &= \sum_{k=0}^{\infty} \left( \frac{\sigma^2}{2} \right)^k \frac{(u^j)^{(2k)}}{k!} \quad \text{(this is a finite sum)}, \\
        &\text{(c)}& \ \frac{\partial \Delta_j (u, \sigma^2)}{\partial u} &= j\Delta_{j-1}(u,\sigma^2) \quad \text{for } j \ge 1, \\
        &\text{(d)}& \ \frac{\partial \Delta_j (u, \sigma^2)}{\partial \sigma^2} &= \frac{j(j-1)}{2} \Delta_{j-2}(u,\sigma^2) \quad \text{for } j \ge 2.
    \end{aligned}
    \end{displaymath}
\end{lemma}

\begin{proof}
    (a): Note that $d \delta_{\sigma^2} / d \xi = - (\xi-u) \delta_{\sigma^2} / {\sigma^2} $. Then for $j \ge 2$ we have
    \begin{displaymath}
    \begin{aligned}
        \Delta_j
        &= \int_{\mathbb{R}} (\xi-u+u) \xi^{j-1} \delta_{\sigma^2} d \xi
        = u \Delta_{j-1} + \int_{\mathbb{R}} (\xi-u) \xi^{j-1} \delta_{\sigma^2} d\xi \\
        &= u \Delta_{j-1} - \sigma^2 \int_{\mathbb{R}} \xi^{j-1} \frac{d \delta_{\sigma^2}}{d\xi}  d\xi
        = u \Delta_{j-1} + (j-1) \sigma^2 \Delta_{j-2}.
    \end{aligned}
    \end{displaymath}
    This, together with  $\Delta_0(u,\sigma^2) = 1$ and $\Delta_1(u,\sigma^2) = u$, indicates that $\Delta_j=\Delta_j(u,\sigma^2)$ is a polynomial of both $u$ and $\sigma^2$.

    (b): This can be proven by induction on $j$. It obviously holds for $\Delta_0=1$ and $\Delta_1=u$. Suppose it is true for $j-1$ and $j$. Then for $j+1$ it follows from (a) that
    \begin{displaymath}
    \begin{aligned}
        \Delta_{j+1}
        &= u\Delta_j + j \sigma^2 \Delta_{j-1}
        = u\sum_{k=0}^{\infty} \left( \frac{\sigma^2}{2} \right)^k \frac{(u^j)^{(2k)}}{k!} + \sigma^2 \sum_{k=0}^{\infty} \left( \frac{\sigma^2}{2} \right)^k \frac{(u^j)^{(2k+1)}}{k!} \\
        &= u^{j+1} + \sum_{k=1}^{\infty} \left( \frac{\sigma^2}{2} \right)^k \frac{u(u^j)^{(2k)}+2k(u^j)^{(2k-1)}}{k!} = \sum_{k=0}^{\infty} \left( \frac{\sigma^2}{2} \right)^k \frac{(u^{j+1})^{(2k)}}{k!}.
    \end{aligned}
    \end{displaymath}
    Hence the proof is complete.

    (c \& d): These two follow immediately from (b):
    \begin{displaymath}
    \begin{aligned}
        \frac{\partial \Delta_{j}}{\partial u}
        &= \sum_{k=0}^{\infty} \left( \frac{\sigma^2}{2} \right)^k \frac{(u^j)^{(2k+1)}}{k!}
        = \sum_{k=0}^{\infty} \left( \frac{\sigma^2}{2} \right)^k \frac{j(u^{j-1})^{(2k)}}{k!}
        = j \Delta_{j-1}; \\
        \frac{\partial \Delta_{j}}{\partial \sigma^2}
        &= \frac{1}{2} \sum_{k=1}^{\infty} k \left( \frac{\sigma^2}{2} \right)^{k-1} \frac{(u^j)^{(2k)}}{k!}
        = \frac{1}{2} \sum_{k=0}^{\infty} \left( \frac{\sigma^2}{2} \right)^{k} \frac{(u^j)^{(2k+2)}}{k!} \\
        &= \frac{1}{2} \sum_{k=0}^{\infty} \left( \frac{\sigma^2}{2} \right)^{k} \frac{j(j-1)(u^{j-2})^{(2k)}}{k!}
        = \frac{j(j-1)}{2} \Delta_{j-2}.
    \end{aligned}
    \end{displaymath}
\end{proof}

\begin{remark} \label{rem:Deltaexplicit}
    \cref{prop:basicDel}(b) is obviously equivalent to
    \begin{equation}
        \Delta_j (u,\sigma^2) = \sum_{k=0}^{[j/2]} \frac{j!}{k!(j-2k)!} \left( \frac{\sigma^2}{2} \right)^k u^{j-2k},
    \end{equation}
    which was established in \cite{MarFox2013}. But the former is more convenient for our later use.
\end{remark}

Inspired by \cref{prop:basicDel}(b), we introduce a family of linear operators $\mathcal{D}_{\vartheta}$, parameterized with $\vartheta \in \mathbb{R}$, acting on the polynomial algebra $\mathbb{R} [u]$. For $f \in \mathbb{R}[u]$, $\mathcal{D}_{\vartheta} f$  is defined as
\begin{equation} \label{eq:opDdefinition}
    \mathcal{D}_{\vartheta} f = \sum_{k=0}^{\infty} \left( \frac{\vartheta}{2} \right)^k \frac{f^{(2k)}}{k!},
\end{equation}
which is a finite sum.
Obviously, $\mathcal{D}_0$ is an identical operator,  $\mathcal{D}_{\sigma^2}f$ is a polynomial of $u$ and $\sigma^2$, and $\mathcal{D}_{\sigma^2} u^j = \Delta_j(u,\sigma^2)$. Further useful properties of $\mathcal{D}_{\vartheta}$ are

\begin{lemma} \label{prop:basicOp}
    \begin{displaymath}
    \begin{aligned}
        &\text{(a)}& \ &\text{(composition)} \quad \mathcal{D}_{\alpha} \circ \mathcal{D}_{\vartheta} = \mathcal{D}_{\alpha+\vartheta}, \\
        &\text{(b)}& \ &\mathcal{D}_{\vartheta} \ \text{is invertible and } \mathcal{D}_{\vartheta}^{-1} = \mathcal{D}_{-\vartheta}, \\
        &\text{(c)}& \ &\frac{\partial}{\partial u} \mathcal{D}_{\vartheta} f(u) = \mathcal{D}_{\vartheta} f'(u), \quad \frac{\partial}{\partial \vartheta} \mathcal{D}_{\vartheta} f(u) = \frac{1}{2} \mathcal{D}_{\vartheta} f''(u), \\
        &\text{(d)}& \ &\mathcal{D}_{\vartheta} (uf) = u\mathcal{D}_{\vartheta}f + \vartheta \mathcal{D}_{\vartheta}f', \\
        &\text{(e)}& \ &\text{If } \mathcal{D}_{\vartheta} f (u_0) = 0, \text{ then } \mathcal{D}_{\vartheta}(uf)|_{u=u_0} = \vartheta \mathcal{D}_{\vartheta} f'(u_0).
    \end{aligned}
\end{displaymath}
\end{lemma}

\begin{proof}
    (a): For the composition, we deduce from the definition that
    \begin{displaymath}
    \begin{aligned}
        (\mathcal{D}_{\alpha} \circ
        &\mathcal{D}_{\vartheta} )f = \sum_{k=0}^{\infty} \left( \frac{\alpha}{2} \right)^k \frac{1}{k!} \sum_{l=0}^{\infty} \left( \frac{\vartheta}{2} \right)^l \frac{f^{(2k+2l)}}{l!} \\
        &= \sum_{p=0}^{\infty} \frac{f^{(2p)}}{p!} \left[ \sum_{l=0}^p \frac{p!}{l!(p-l)!} \left( \frac{\alpha}{2} \right)^{p-l} \left( \frac{\vartheta}{2} \right)^l \right] 
        = \sum_{p=0}^{\infty} \frac{f^{(2p)}}{p!} \left( \frac{\alpha+\vartheta}{2} \right)^p = \mathcal{D}_{\alpha+\vartheta} f.
    \end{aligned}
    \end{displaymath}

    (b) follows immediately from (a) and $\mathcal{D}_0 = id$.

    For (c), the first one is obvious, while the second can be shown as \cref{prop:basicDel}(d).

    (d): By using $(uf)^{(2k)} = u f^{(2k)} + 2k f^{(2k-1)}$, this can be proved as \cref{prop:basicDel}(b). Then (e) follows immediately from (d).
\end{proof}

\subsection{Hyperbolicity of EQMOM}
\label{subsec:hypeqmom}

In this section we prove that the EQMOM-derived moment system \cref{eq:EQMOMboltz} for the 1-D Boltzmann equation is strictly hyperbolic, which will be shown to be sufficient for the structural stability condition (ii). The conclusion can be stated as

\begin{theorem} \label{thm:hypeqmom}
    For $M \in \mathcal{M}(\Omega_W)$, the $(2N+1) \times (2N+1)$ coefficient matrix $A=A(M)$ in \cref{eq:EQMOMA} has $(2N+1)$ distinct real eigenvalues. Namely, the EQMOM-derived moment system \cref{eq:EQMOMboltz} is strictly hyperbolic.
\end{theorem}

We should mention that this theorem was already established in \cite{Chalons2017} for $N=2$ (the two-node system) but the proof does not seem to work for $N > 2$.

Our proof of this theorem needs some preparations. First of all, the characteristic polynomial of $A$ in \cref{eq:EQMOMA} reads as
\begin{equation} \label{eq:Acharpeqmom}
    c(u;W) = u^{2N+1} - a_{2N}u^{2N} - \cdots - a_1 u - a_0.
\end{equation}
Here the coefficient $a_j =a_j(W)= \frac{\partial \bar{M}_{2N+1}}{\partial M_j}$ ($j=0,1,\dots,2N$), with $\bar{M}_{2N+1}$ defined in \cref{eq:EQMOMM_2Np1}, is a function of $W$.
To show that $c(u;W)$, as a polynomial of $u$, has $(2N+1)$ distinct real roots for $M \in \mathcal{M}(\Omega_W)$, we introduce an auxiliary function
\begin{equation} \label{eq:gdefine}
    g(u;W) = \mathcal{D}_{\sigma^2} c(u;W) = \sum_{k=0}^{\infty} \left( \frac{\sigma^2}{2} \right)^k \frac{\partial_u^{2k} c(u;W)}{k!}.
\end{equation}
By \cref{prop:basicOp}(b), we have
\begin{equation} \label{eq:cbyg}
    c(u;W) = \mathcal{D}_{-\sigma^2} g(u;W) = \sum_{k=0}^{\infty} \left( - \frac{\sigma^2}{2} \right)^k \frac{\partial_u^{2k} g(u;W)}{k!}.
\end{equation}

Set $a_{2N+1}=-1$. Then $c(u;W)$ can be rewritten as $-\sum_{j=0}^{2N+1} a_j u^{j}$ and from the linearity of $\mathcal{D}_{\sigma^2}$ it follows that
\begin{equation} \label{eq:ggg}
    g(u;W) = -\sum_{j=0}^{2N+1} a_j \mathcal{D}_{\sigma^2} u^{j} = -\sum_{j=0}^{2N+1} a_j \Delta_j(u, \sigma^2).
\end{equation}
Moreover, from \cref{eq:Acharpeqmom,eq:gdefine} we see that $g(u;W)$ is a $u$-polynomial of degree $(2N+1)$:
\begin{equation} \label{eq:cgrelate1}
    g(u;W) = -\sum_{j=0}^{2N+1} g_j u^j
\end{equation}
with $g_{2N+1} = a_{2N+1}= -1$. Further relations between the coefficients of $g(u;W)$ and $c(u;W)$ are
\begin{equation} \label{eq:ajbygj}
    a_j = \sum_{k=0}^{N-[j/2]} g_{j+2k} \frac{(j+2k)!}{j!k!} \left( - \frac{ \sigma^2}{2} \right)^k, \quad j=0,1,\dots, 2N+1.
\end{equation}
This can be shown as
\begin{displaymath}
\begin{aligned}
        c(u;W)
        &= \sum_{k=0}^N \frac{\partial_u^{2k} g(u;W)}{k!} \left( -\frac{\sigma^2}{2} \right)^k
        = - \sum_{k=0}^N \sum_{j=0}^{2N+1-2k} \frac{(j+2k)!}{j!k!} \left( - \frac{\sigma^2}{2} \right)^k g_{j+2k} u^j \\
        &= - \sum_{j=0}^{2N+1} \left[ \sum_{k=0}^{N-[j/2]} \frac{(j+2k)!}{j!k!} \left( -\frac{\sigma^2}{2} \right)^k g_{j+2k} \right] u^j.
\end{aligned}
\end{displaymath}

Furthermore, $g=g(u)=g(u;W)$ has the following elegant expression.

\begin{lemma} \label{prop:groots}
    \begin{displaymath}
        g(u;W) = (u-u_1)^2 \cdots (u-u_N)^2 (u-\tilde{U}),
    \end{displaymath}
    where $u_1,\dots,u_N$ are the nodes solved from \cref{eq:EQMOMmo}, and
    \begin{displaymath}
        \tilde{U} = \tilde{U} (W) =  \frac{\sum_{i=1}^N w_i u_i \prod_{1\le j\le N, j\ne i} (u_j-u_i)^2}{\sum_{i=1}^N w_i \prod_{1\le j\le N, j \ne i} (u_j-u_i)^2}
    \end{displaymath}
    for $W \in \Omega_W^{open}$.
\end{lemma}

\begin{remark} \label{rem:geq}
    This lemma shows that for $W \in \Omega_{W}^{open}$, $\tilde{U}$ is a convex combination of the $u_i$'s. Moreover, for $W = (w_1, U, w_2, U, \dots, w_N, U, \sigma^2) \in \Omega_{W}^{eq}$ and any sequence $\{W_k\} \subset \Omega_{W}^{open}$ approaching $W$, $\tilde{U}(W_k)$ converges to $U$.
    Because of this, for $W \in \Omega_W^{eq}$ we define $\tilde{U}(W) = U$ ($=M_1/M_0$) and thereby $g(u;W) = (u-U)^{2N+1}$.
\end{remark}

\begin{proof}
    By \cref{prop:basicDel}(c\&d), the Jacobian matrix of the map $M = \mathcal{M}(W)$ defined in \cref{eq:EQMOMmo} is
    \begin{equation} \label{eq:jacoeqmom}
        \begin{bmatrix}
            \Delta_0(u_1) & 0 & \cdots & \Delta_0(u_N) & 0 & 0 \\
            \Delta_1(u_1) & w_1 \Delta_0(u_1) & \cdots & \Delta_1(u_N) & w_N \Delta_1(u_N) & 0 \\
            \vdots & \vdots & & \vdots & \vdots & \vdots \\
            \Delta_j(u_1) & jw_1 \Delta_{j-1}(u_1) & \cdots & \Delta_j(u_N) & j w_N \Delta_{j-1}(u_N) & \binom{j}{2} M_{j-2} \\
            \vdots & \vdots & & \vdots & \vdots & \vdots \\
            \Delta_{2N}(u_1) & 2N w_1 \Delta_{2N-1}(u_1) & \cdots & \Delta_{2N}(u_N) & 2N w_N \Delta_{2N-1}(u_N) & \binom{2N}{2} M_{2N-2}
        \end{bmatrix}.
    \end{equation}
    Note that the dependence of $\Delta_j$ on $\sigma^2$ has been omitted here for clarity. Moreover, from \cref{eq:EQMOMM_2Np1}, $\partial \bar{M}_{2N+1} / \partial W$ reads as
    \begin{displaymath}
        \left( \Delta_{2N+1}(u_1), (2N+1)w_1 \Delta_{2N}(u_1), \dots, \Delta_{2N+1}(u_N), (2N+1)w_N \Delta_{2N}(u_N), \binom{2N+1}{2} M_{2N-1} \right).
    \end{displaymath}
    Then from the simple relation
    \begin{displaymath}
        (a_0,a_1,...,a_{2N}) \frac{\partial \mathcal{M}}{\partial W} = \left( \frac{\partial \bar{M}_{2N+1}}{\partial M} \right) \left( \frac{\partial \mathcal{M}}{\partial W} \right) = \frac{\partial \bar{M}_{2N+1}}{\partial W}
    \end{displaymath}
    we obtain
    \begin{subequations}
    \begin{align}
        a_0 \Delta_0 (u_k) + a_1 \Delta_1 (u_k) + \dots + a_{2N} \Delta_{2N} (u_k) &= \Delta_{2N+1} (u_k), \label{eq:aMWeqmom1} \\
        w_k a_1 \Delta_0 (u_k) + \dots + 2N w_k a_{2N} \Delta_{2N-1} (u_k) &= (2N+1) w_k \Delta_{2N} (u_k), \label{eq:aMWeqmom2} \\
        \binom{2}{2} a_2 M_0 + \dots + \binom{2N}{2} a_{2N} M_{2N-2} &= \binom{2N+1}{2} M_{2N-1}. \label{eq:aMWeqmom3}
    \end{align}
    \end{subequations}
    \cref{eq:aMWeqmom1,eq:aMWeqmom2}, together with \cref{eq:ggg} and \cref{prop:basicDel}(c), imply that $g(u_k)=0$ and $\left. \frac{dg(u)}{du} \right|_{u=u_k} =0$ for $k=1,\dots,N$.
    Thus we see the expected expression of $g(u;W)$ from \cref{eq:cgrelate1} with $\tilde{U}$ to be determined.

    Next, we use \cref{eq:aMWeqmom3} to determine $\tilde{U}$.
    Recall that $a_{2N+1}=-1$. We use \cref{eq:EQMOMmo} to rewrite \cref{eq:aMWeqmom3} as
    \begin{equation} \label{eq:ajsrel}
        0=\sum_{j=2}^{2N+1} \binom{j}{2} a_j M_{j-2} = \sum_{i=1}^N w_i \sum_{j=2}^{2N+1} \binom{j}{2} a_j \Delta_{j-2} (u_i,\sigma^2).
    \end{equation}
    On the other hand, we deduce from \cref{eq:ggg,eq:cgrelate1} that
    \begin{displaymath}
        -\frac{1}{2} g''(u) = \sum_{j=2}^{2N+1} \binom{j}{2} g_j u^{j-2} = \sum_{j=2}^{2N+1} \binom{j}{2} a_j \Delta_{j-2}(u,\sigma^2).
    \end{displaymath}
    Then we see from \cref{eq:ajsrel} that
    \begin{equation} \label{eq:ajsrel1}
        \sum_{i=1}^{N} w_i \sum_{j=2}^{2N+1} \binom{j}{2} g_j u_i^{j-2} = 0.
    \end{equation}

    Now we define $\tilde{g}(u) = (u-u_1)^2 \cdots (u-u_N)^2 = -\sum_{j=0}^{2N} \tilde{g}_j u^{j}$ with $\tilde{g}_{2N} = -1$. Then $g(u) = (u-\tilde{U}) \tilde{g}(u)$ and the coefficients are related with
    \begin{displaymath}
        g_j = \tilde{g}_{j-1} - \tilde{U} \tilde{g}_j
    \end{displaymath}
    for $0 \le j \le 2N+1$ ($\tilde{g}_{-1} = \tilde{g}_{2N+1} = 0$). Substituting this relation into \cref{eq:ajsrel1}, we obtain
    \begin{displaymath}
        \left[ \sum_{j=2}^{2N} \binom{j}{2} \tilde{g}_j M_{j-2}^{*} \right] \tilde{U} = \sum_{j=2}^{2N+1} \binom{j}{2} \tilde{g}_{j-1} M_{j-2}^{*},
    \end{displaymath}
    where $M_j^{*} = \sum_{i=1}^{N} w_i u_i^{j}$. It remains to show
    \begin{displaymath}
    \begin{aligned}
        \sum_{j=2}^{2N} \binom{j}{2} \tilde{g}_j M_{j-2}^{*} &= - \sum_{i=1}^N w_i \prod_{1\le k\le N,k \ne i} (u_k-u_i)^2, \\
        \sum_{j=2}^{2N+1} \binom{j}{2} \tilde{g}_{j-1} M_{j-2}^{*} &= - \sum_{i=1}^N w_i u_i \prod_{1\le k\le N, k \ne i} (u_k-u_i)^2.
    \end{aligned}
    \end{displaymath}
    These two follow from the obvious relations
    \begin{displaymath}
    \begin{aligned}
        - \sum_{j=2}^{2N} \binom{j}{2} \tilde{g}_j u_i^{j-2} &= \frac{1}{2} \tilde{g}''(u_i) = \prod_{1\le k\le N, k \ne i} (u_k - u_i)^2, \label{eq:gtiltsecder} \\
        - \sum_{j=2}^{2N+1} \binom{j}{2} \tilde{g}_{j-1} u_i^{j-2} &= \frac{1}{2} \left( u \tilde{g} \right)'' (u_i) = u_i \prod_{1\le k\le N, k \ne i} (u_k-u_i)^2 \label{eq:ugtiltsecder}
    \end{aligned}
    \end{displaymath}
    for any $1 \le i \le N$.
    This completes the proof.
\end{proof}

\begin{remark} \label{rem:greqgorder}
    \cref{prop:groots} indicates that the coefficients $g_j$ of $g(u;W)$ in \cref{eq:cgrelate1} are independent of $\sigma^2$. From \cref{eq:cbyg} we see that $c(u;W)$ is a bivariate polynomial of $u$ and $\sigma^2$, the coefficients $a_j$ of $c(u;W)$ are polynomials of $\sigma^2$, and  $c(u;W) = g(u)$ for $\sigma^2=0$.
    Furthermore, the $j$th derivative $c^{(j)}(u;W)$ of $c(u;W)$ with respect to $u$ can be viewed as a perturbation of $g^{(j)}(u)$ with the single parameter $\sigma^2 \ge 0$ for $0 \le j \le 2N+1$.
\end{remark}

By \cref{prop:groots}, $g(u)$ has $(2N+1)$ ($=$ the degree of $g$) real roots (including multiplicity). This fact can be further generalized as follows.
\begin{lemma} \label{prop:gderroots}
    For any $0 \le j \le 2N$, $g^{(j)}(u)$ has $(2N+1-j)$ real roots (including multiplicity). Hence any local minimum (maximum) value of $g^{(j)}(u)$ is non-positive (non-negative).
\end{lemma}

\begin{proof}
    We prove by induction on $j$. As discussed above, the conclusion holds for $j=0$. Namely, $g$ has $(2N+1)$ roots. Suppose it holds for $0,\dots,j$. Then we have $g^{(j)}(u) = C (u-\tilde{u}_1)^{k_1} \cdots (u-\tilde{u}_m)^{k_m}$, where $m\ge 1$, $\tilde{u}_1 < \tilde{u}_2 < \dots < \tilde{u}_m$, $k_i \ge 1$ and $k_1 + \dots + k_m = 2N+1-j$.
    Thus $(u-\tilde{u}_i)^{k_i-1}$ is a factor of $g^{(j+1)}(u)$ for any $1\le i \le m$.
    Besides, Rolle's theorem implies the existence of at least one root of $g^{(j+1)}(u)$ in each open interval $(\tilde{u}_i,\tilde{u}_{i+1})$ for $1\le i \le m-1$. Therefore, the number of roots of $g^{(j+1)}$ is no less than
\begin{displaymath}
    (k_1-1)+\cdots+(k_m-1) + (m-1) = (k_1+\cdots+k_m)-1 = 2N-j.
\end{displaymath}
Since $g^{(j+1)}$ is of degree $(2N-j)$, it must have $(2N-j)$ roots (including multiplicity). This also indicates that $g^{(j)}$ has only one extreme point in each open interval above. Hence any local minimum (maximum) value of $g^{(j)}(u)$ is non-positive (non-negative).
\end{proof}

With the preparations above, we are in a position to prove \cref{thm:hypeqmom}.

\begin{proof}[Proof of \cref{thm:hypeqmom}]
    We will prove the following stronger statement: for $0\le j\le 2N+1$, $c^{(2N+1-j)}(u;W)$ has $j$ distinct roots for any $W =(w_1,u_1,\dots, w_N,u_N,\sigma^2)$ $\in \Omega_W$ with $\sigma^2>0$. This will be done with induction on $j$.
    For $j=0, \ 1$, the statement is obvious because $c^{(2N+1)}(u;W) = (2N+1)!$ and $c^{(2N)}(u;W)$ is of degree 1.

    Suppose the conclusion holds for $j \le k (\le 2N+1)$.
    From \cref{rem:greqgorder} we know that $c^{(2N+1-k)}(u;W)$ is a bivariate polynomial of $u$ and $\sigma^2$ on $\mathbb{R} \times [0,\infty)$.
    Denote $u^{*}(\sigma^2) \in \mathbb{R}$ to be one root of $c^{(2N+1-k)}(u;W)$. Thus $u^{*}(\sigma^2)$ is an extreme point of $c^{(2N-k)}(u;W)$ and $u^{*}(0)$ is a root of $g^{(2N+1-k)}(u)$.
    Moreover, $u^{*}(\sigma^2)$ is continuous on $\sigma^2 \in [0,\infty)$ and differentiable on $(0,\infty)$ because the roots are distinct \cite{Kato1980}.

    Next we consider the extreme values of $c^{(2N-k)}(u;W)$ at $u=u^*(\sigma^2)$.
    Since $c^{(2N-k)}(u;W)$ is a polynomial of $u$ and $\sigma^2$, the composite $h_{k}(\sigma^2):=c^{(2N-k)}(u^*(\sigma^2);W)$ is continuous on $[0,\infty)$ and differentiable on $(0,\infty)$. And $h_k(0)$ is the extreme value of $g^{(2N-k)}(u)$.
    According to \cref{prop:gderroots}, $h_k(0) \ge 0$ if it is a local maximum and $h_k(0) \le 0$ if it is a local minimum.
    For $\sigma^2>0$, because $c^{(2N+1-k)}(u^*(\sigma^2);W)=0$, the derivative of $h_k(\sigma^2)$ reads as
    \begin{displaymath}
    \begin{aligned}
        \frac{\partial h_k(\sigma^2)}{\partial \sigma^2}
        &= \frac{\partial}{\partial \sigma^2} c^{(2N-k)}(u^*(\sigma^2);W)
        = \frac{\partial}{\partial \sigma^2} \sum_{l=0}^{\infty} \left( - \frac{\sigma^2}{2} \right)^l \frac{g^{(2l+2N-k)}(u^*(\sigma^2))}{l!} \\
        &= -\frac{1}{2} \sum_{l=1}^{\infty} \left( - \frac{\sigma^2}{2} \right)^{l-1} \frac{g^{(2l+2N-k)}(u^*(\sigma^2))}{(l-1)!} + c^{(2N+1-k)}(u^*(\sigma^2);W)
        \frac{\partial u^*(\sigma^2)}{\partial \sigma^2} \\
        &= -\frac{1}{2} c^{(2N+2-k)}(u^*(\sigma^2);W).
    \end{aligned}
    \end{displaymath}
    Thus, if $c^{(2N-k+2)}(u^*(\sigma^2);W)<0$ (that is, $u^*(\sigma^2)$ is a local maximum point of $c^{(2N-k)}(u;W)$), then the local maximum value $h_k(\sigma^2)$ strictly increases on $\sigma^2 \in (0,\infty)$.
    Since $h_k(0) \ge 0$ and $h_k(\sigma^2)$ is continuous at $\sigma^2=0$, we conclude that $h_k(\sigma^2)>0$ for all $\sigma^2>0$.
    Similarly, if $c^{(2N-k+2)}(u^*(\sigma^2);W)>0$, we have $h_k(\sigma^2)<0$ for all $\sigma^2>0$.

    In summary, the above arguments show that each local maximum value of the ($k+1$)th oder polynomial $c^{(2N-k)}(u;W)$ is positive and each local minimum value is negative. On the other hand, by the induction assumption $c^{(2N+1-k)}(u;W)$ has $k$ distinct real roots, which are naturally extreme points of $c^{(2N-k)}(u;W)$ for $\sigma^2>0$.  Therefore, $c^{(2N-k)}(u;W)$ has $(k-1)$ distinct real roots among the extreme points. Moreover, the induction assumption implies that $c^{(2N-k)}(u;W)$ has one root larger and another one less than all the extreme points.  Thus, for each $\sigma^2>0$, $c^{(2N-k)}(u;W)$ has $(k+1)$ distinct real roots. By the induction principle this completes the proof.
\end{proof}

\begin{remark} \label{rem:hyp2}
    By \cref{thm:hypeqmom}, the coefficient matrix $A=A(M)$ of the 1-D moment system \cref{eq:EQMOMboltz} has $n=(2N+1)$ distinct real eigenvalues $\lambda_i$ ($1\le i \le n$) for $\sigma^2>0$. Denote by $r_i$ the corresponding left eigenvectors. Set $L = (r_1^{T},\dots,r_n^{T})^{T}$.
    It is clear that $A_0(M) = L^{T} \Lambda L$ with $\Lambda$ an arbitrary positive diagonal matrix is a symmetrizer in the structural stability condition (ii).
    As a matter of fact,  it is straightforward to show that such a symmetrizer can only be of the form $L^{T} \Lambda L$ .
\end{remark}

\subsection{Equilibrium state}
\label{subsec:eqmomequistate}

As stated in \Cref{subsec:stability}, (i) and (iii) of the structural stability condition should be examined on the equilibrium manifold $\mathcal{E}$ where $S(M(W))=0$.
In this section we determine the equilibrium manifold.

For the BGK model, $S(M(W))=0$ is equivalent to
\begin{equation} \label{eq:equieq}
    \sum_{i=1}^{N} w_i \Delta_j (u_i, \sigma^2) = M_j = \rho \Delta_j (U, \theta) \quad \text{for } j=0,\dots,2N
\end{equation}
(see \Cref{subsubsec:EQMOM}). Thus, the equilibrium state $W=(w_1, u_1, \dots, w_N, u_N, \sigma^2)^T$ is determined by the three macroscopic parameters $\rho$, $U$ and $\theta$. And we need to find $W$ from \cref{eq:equieq} for $1 \le i \le N$.

For this purpose, we recall from \Cref{subsec:preli} that $\Delta_0(u,\sigma^2)=1$, $\Delta_1(u,\sigma^2)=u$ and $\Delta_2(u,\sigma^2)=u^2+\sigma^2$. Thus, for $j=0,1,2$, \cref{eq:equieq} is just
\begin{displaymath}
    \sum_{i=1}^{N} w_i = \rho, \quad \sum_{i=1}^{N} w_i u_i = \rho U, \quad \sum_{i=1}^{N} w_i u_i^2 = \rho U^2 + \rho (\theta - \sigma^2).
\end{displaymath}
Then we deduce from the inequality $\left( \sum_{i=1}^{N} w_i \right) \left( \sum_{i=1}^{N} w_i u_i^2 \right) \ge \left( \sum_{i=1}^{N} w_i u_i \right)^2$
that
\begin{equation} \label{eq:zetaprop}
    \sigma^2 \le \theta \quad \text{and } \sigma^2=\theta \text{ if and only if all the $u_i$'s are equal}.
\end{equation}

For further discussions, we need the following fact.
\begin{proposition} \label{lem:Mjstar}
    \begin{displaymath}
        M_j^{*}: = \sum_{i=1}^{N} w_i u_i^{j}= \sum_{k=0}^{[j/2]} \frac{j!}{k!(j-2k)!} \left( - \frac{\sigma^2}{2} \right)^k M_{j-2k}.
    \end{displaymath}
\end{proposition}

\begin{proof}
    Recall that $\Delta_j(u,\sigma^2) = \mathcal{D}_{\sigma^2} u^j$. From \cref{prop:basicOp}(b) and \cref{prop:basicDel}(c) we deduce that
    \begin{displaymath}
        u^j = \sum_{k=0}^{\infty} \frac{\partial_u^{2k} \Delta_j(u,\sigma^2)}{k!} \left( - \frac{\sigma^2}{2} \right)^k = \sum_{k=0}^{[j/2]} \frac{j!}{k!(j-2k)!} \left( - \frac{\sigma^2}{2} \right)^k \Delta_{j-2k}(u,\sigma^2).
    \end{displaymath}
    Then taking the weighted summation $\sum_{i=1}^N w_i$ and using \cref{eq:EQMOMmo} give the proposition.
\end{proof}

Next we define $\zeta^2 = \theta - \sigma^2 \ge 0$ and show that \cref{eq:equieq} is equivalent to
\begin{equation} \label{eq:equieqred}
    \sum_{i=1}^{N} w_i u_i^j = \rho \Delta_j (U, \zeta^2) \quad \text{for } j=0,\dots,2N.
\end{equation}
Indeed, if \cref{eq:equieq} holds (i.e. $M(W) \in \mathcal{E}$), the last proposition implies that
\begin{displaymath}
\begin{aligned}
    \sum_i w_i u_i^j
    = M_j^* &= \sum_{k=0}^{[j/2]} \left( -\frac{\sigma^2}{2} \right)^k \frac{j!}{k!(j-2k)!} \left[ \rho \Delta_{j-2k} (U,\theta) \right] \\
    &= \rho \sum_{k=0}^{\infty} \left( -\frac{\sigma^2}{2} \right)^k \left. \frac{ \partial_u^{2k} \Delta_j (u,\theta) }{k!} \right|_{u=U}
    = \rho \mathcal{D}_{-\sigma^2} \mathcal{D}_{\theta} U^j = \rho \mathcal{D}_{\theta-\sigma^2} U^j.
\end{aligned}
\end{displaymath}
Here the expression $\mathcal{D}_{\vartheta} f(U)$ denotes $\mathcal{D}_{\vartheta} f(u)|_{u=U}$ for arbitrary polynomial $f$ and the last step is due to \cref{prop:basicOp}(a). This is just \cref{eq:equieqred}.
The deduction of \cref{eq:equieq} from \cref{eq:equieqred} is similar.

Now we are in a position to state the central result of this section.
\begin{theorem} \label{thm:equisolu}
    The equilibrium state belongs to $\Omega_W^{eq}$, that is,
    \begin{displaymath}
        u_1 = \dots = u_N = U, \ \sigma^2 = \theta, \ \text{and} \ \sum_{i=1}^{N} w_i = \rho.
    \end{displaymath}
    Hence, at equilibrium $\bar{M}_{2N+1} = \rho \Delta_{2N+1} (U, \theta)$.
\end{theorem}

\begin{proof}
    Thanks to \cref{eq:zetaprop}, it suffices to show that $\zeta^2 := \theta - \sigma^2 = 0$. Otherwise, the $u_i$'s must take $N'$ different values ($1<N'\le N$).
    Then, by redefining $w_i$, the summation $\sum_{i=1}^N w_i u_i^j$ in the left-hand side of \cref{eq:equieqred} is reduced to $\sum_{k=1}^{N'} w_k u_k^j$ where the $u_k$'s are distinct ($1 \le k \le N'$).
    Thus, we may as well assume that all the $u_i$'s are distinct and $\zeta^2>0$.
    Then we will derive a contradiction in three steps, where the abbreviation $$\mathcal{D}_{\vartheta} f(U) \equiv \mathcal{D}_{\vartheta} f(u)|_{u=U}$$
    will be frequently used.

    \textbf{Step I}. Because $\Delta_j (u,\zeta^2) = \mathcal{D}_{\zeta^2} u^j$, the first $N$ equations ($j=0,\dots,N-1$) in \cref{eq:equieqred} can be rewritten as a system of linear algebraic equations:
    \begin{displaymath}
        \begin{bmatrix}
            1& 1 & \cdots & 1 \\
            u_1 & u_2 & \cdots & u_N \\
            \vdots & \vdots & & \vdots \\
            u_1^{N-1} & u_2^{N-1} & \cdots & u_N^{N-1}
        \end{bmatrix}
        \begin{bmatrix}
            w_1 \\ w_2 \\ \vdots \\ w_N
        \end{bmatrix}
        = \rho
        \begin{bmatrix}
            \mathcal{D}_{\zeta^2} (1) \\
            \mathcal{D}_{\zeta^2} (U^1) \\
            \vdots \\
            \mathcal{D}_{\zeta^2} (U^{N-1})
        \end{bmatrix}.
    \end{displaymath}
    Since all the $u_i$'s are distinct, this gives a unique $(w_1,\dots,w_N)$ in terms of $(u_1,\dots,u_N)$, $\rho$, $U$ and $\zeta^2$. We claim that for $1\le i \le N$,
    \begin{equation} \label{eq:wbyueq}
        w_i \prod_{1\le k \le N, k \ne i} (u_i-u_k) = \rho \mathcal{D}_{\zeta^2} \left( \prod_{1 \le k \le N, k \ne i} (U-u_k) \right).
    \end{equation}
    To see this, we use the uniqueness and only need to show that the $w_i$'s solve the system of equations above.
    Indeed, thanks to the Lagrange interpolating polynomial
    \begin{equation*} \label{eq:lip1}
        \sum_{i=1}^N \frac{ \prod_{1 \le k \le N, k \ne i} (u-u_k) }{ \prod_{1 \le k \le N, k \ne i} (u_i-u_k) } u_i^j = u^j \quad \text{for } 0 \le j \le N-1
    \end{equation*}
    and the linearity of the operator $\mathcal{D}_{\zeta^2}, $\cref{eq:wbyueq} implies that for $0\le j \le N-1$,
    \begin{displaymath}
        \sum_{i=1}^N w_i u_i^j
        = \rho \mathcal{D}_{\zeta^2} \left( \sum_{i=1}^N \frac{ \prod_{1 \le k \le N, k \ne i} (U-u_k) }{ \prod_{1 \le k \le N, k \ne i} (u_i-u_k) } u_i^j \right) = \rho \mathcal{D}_{\zeta^2} (U^j).
    \end{displaymath}
    Namely, the $w_i$'s defined in \cref{eq:wbyueq} solve the system of linear algebraic equations above.

    \textbf{Step II}. With the $w_i$'s defined in \cref{eq:wbyueq}, we turn to the next $N$ equations ($j=N,\dots,2N-1$) in \cref{eq:equieqred} to solve $u_i$:
    \begin{displaymath}
        \begin{bmatrix}
            1& 1 & \cdots & 1 \\
            u_1 & u_2 & \cdots & u_N \\
            \vdots & \vdots & & \vdots \\
            u_1^{N-1} & u_2^{N-1} & \cdots & u_N^{N-1}
        \end{bmatrix}
        \begin{bmatrix}
            w_1 u_1^N \\ w_2 u_2^N \\ \vdots \\ w_N u_N^N
        \end{bmatrix}
        = \rho
        \begin{bmatrix}
            \mathcal{D}_{\zeta^2} (U^N) \\
            \mathcal{D}_{\zeta^2} (U^{N+1}) \\
            \vdots \\
            \mathcal{D}_{\zeta^2} (U^{2N-1})
        \end{bmatrix}.
    \end{displaymath}
    Again, the solution $w_iu_i^N$ is unique. As in \textbf{Step I}, we can show that
    \begin{equation} \label{eq:wuNeq}
        w_i u_i^N \prod_{1\le k \le N, k \ne i} (u_i-u_k) = \rho \mathcal{D}
        _{\zeta^2} \left( U^N \prod_{1\le k\le N, k\ne i} (U-u_k) \right)
    \end{equation}
    for $1\le i \le N$.

    Substituting \cref{eq:wbyueq} into \cref{eq:wuNeq}, we obtain
    \begin{equation} \label{eq:usolve1}
        u_i^N \mathcal{D}_{\zeta^2} \left( \prod_{1\le k \le N, k\ne i} (U-u_k) \right) = \mathcal{D}_{\zeta^2} \left( U^N \prod_{1\le k \le N, k\ne i} (U-u_k) \right)
    \end{equation}
    for $1 \le i \le N$. By the linearity of $\mathcal{D}_{\zeta^2}$, \cref{eq:usolve1} is equivalent to
    \begin{displaymath}
        \mathcal{D}_{\zeta^2} \left( (U^{N-1} + u_i U^{N-2} + \cdots + u_i^{N-1}) \prod_{k=1}^N (U-u_k) \right) = 0,
    \end{displaymath}
    which can be rewritten as
    \begin{displaymath}
        \begin{bmatrix}
            1 & u_1 & \cdots & u_1^{N-1} \\
            \vdots & \vdots & & \vdots \\
            1 & u_N & \cdots & u_N^{N-1}
        \end{bmatrix}
        \begin{bmatrix}
            \mathcal{D}_{\zeta^2} (U^{N-1} F) \\
            \vdots \\
            \mathcal{D}_{\zeta^2} (F)
        \end{bmatrix} = 0
    \end{displaymath}
    with $F=F(U) = \prod_{k=1}^N (U-u_k)$.
    Since all the $u_i$'s are distinct, this says
    \begin{equation} \label{eq:dFeq0}
        \mathcal{D}_{\zeta^2}(F) = \mathcal{D}_{\zeta^2}(UF) = \cdots = \mathcal{D}_{\zeta^2}(U^{N-1}F) = 0.
    \end{equation}

    Having this, in \cref{prop:basicOp}(e) we take $u_0 = U$ and $f(u) = u^j F(u)$ ($0 \le j \le N-2$) and deduce from \cref{eq:dFeq0} that
    $$0 = \mathcal{D}_{\zeta^2} (U^{j+1}F) = \zeta^2 \mathcal{D}_{\zeta^2} \left( (U^jF)' \right) = \zeta^2 \mathcal{D}_{\zeta^2} \left( U^j F' \right).$$
    Hence $\mathcal{D}_{\zeta^2}(U^j F') = 0$ for $0 \le j \le N-2$. This procedure can be repeated for the derivative of $F'(u)$ to yield $\mathcal{D}_{\zeta^2} (U^j F'') = 0$ for $0 \le j \le N-3$. Moreover, we have
    \begin{equation} \label{eq:dFjeq0}
        \mathcal{D}_{\zeta^2} \left( U^j F^{(k)} \right) = 0
    \end{equation}
    for $0 \le k \le N-1$ and $0 \le j \le N-1-k$.

    \textbf{Step III}. In this step, we use \cref{eq:wuNeq} and the Lagrange interpolating polynomial
    \begin{equation*} \label{eq:interpj}
            \sum_{i=1}^N \frac{\prod_{1\le k\le N, k\ne i} (u-u_k)}{\prod_{1\le k\le N, k\ne i} (u_i-u_k)} u_i^N
            = u^N - \prod_{k=1}^N (u-u_k)
    \end{equation*}
    to deduce that
    \begin{displaymath}
    \begin{aligned}
        \sum_{i=1}^N w_i u_i^{2N}
        &= \rho \mathcal{D}_{\zeta^2} \left( U^N \sum_{i=1}^N \frac{\prod_{1\le k\le N, k\ne i} (U-u_k)}{\prod_{1\le k\le N, k\ne i} (u_i-u_k)} u_i^N \right) = \rho \mathcal{D}_{\zeta^2} \left( U^N(U^N-F) \right).
    \end{aligned}
    \end{displaymath}
    Thus, the last equation in \cref{eq:equieqred} is equivalent to $\mathcal{D}_{\zeta^2} (U^N (U^N-F)) = \mathcal{D}_{\zeta^2} (U^{2N})$ or
    \begin{equation*}
        \mathcal{D}_{\zeta^2} \left( U^N F \right) = 0.
    \end{equation*}
    Then we use \cref{prop:basicOp}(d) and \cref{eq:dFjeq0} to see that
    $$0 = \mathcal{D}_{\zeta^2} (U^N F) = \zeta^2 \mathcal{D}_{\zeta^2} (U^{N-1} F') =
    \cdots = \zeta^{2N} \mathcal{D}_{\zeta^2} (F^{(N)}) = N! \cdot \zeta^{2N},$$
    which implies that $\zeta^2=0$.
    This contradicts the assumption that all the $u_i$'s are distinct and $\zeta^2>0$.
    Hence the proof is complete.
\end{proof}

\subsection{BGK model}
\label{subsec:coupAseqmom}

In this subsection we show that the EQMOM moment system \cref{eq:EQMOMboltz} with BGK source term
$$S(M)= \rho \left(\Delta_0 (U,\theta),\dots, \Delta_{2N} (U,\theta) \right)^{T} - M$$
satisfies the structural stability condition (i)--(iii). Indeed, (ii) has been verified in \cref{rem:hyp2}.

To see Condition (i), we compute the Jacobian matrix of $S=S(M)$. Notice that the first three components of $S$ vanish identically and $\rho, U,\theta$ depend only on $M_0$, $M_1$ and $M_2$.
Then the Jacobian matrix can be written as
\begin{equation} \label{eq:sourcejac}
    S_M (M) := \frac{\partial S}{\partial M} =
    \begin{bmatrix}
        0_{3 \times 3} &  \\
        \hat{S}_M & -I_{2N-2}
    \end{bmatrix},
\end{equation}
where $\hat{S}_M$ is a $(2N-2)\times 3$ matrix with
\begin{equation} \label{eq:chidef}
    \left( \hat{S}_M \right)_{i-2, \ j+1} = \chi_i^j: = \partial (\rho \Delta_i(U,\theta)) / \partial M_j
\end{equation}
for $3 \le i \le 2N$ and $j=0,1,2$.
Now we take
\begin{equation} \label{eq:sourcediag}
    P=\begin{bmatrix} I_3 & \\ - \hat{S}_M & I_{2N-2} \end{bmatrix}
\end{equation}
and see that $P S_M = \begin{bmatrix} 0_{3\times 3}& \\ & -I_{2N-2} \end{bmatrix} P$, which justifies Condition (i).

The rest of this subsection is to show Condition (iii).
To this end, we need to choose the symmetrizer $A_0 = A_0(M)$.
As pointed out in \cref{rem:hyp2}, such a symmetrizer $A_0$ can only be of the form $L^{T} \Lambda L$ with $\Lambda$ a diagonal positive definite matrix to be determined.

Firstly, we specify the matrix $L = (r_1^{T},\dots,r_{2N+1}^{T})^{T}$ with $r_i$ a left eigenvector of the coefficient matrix $A=A(M)$ corresponding to the eigenvalues $\lambda_i$ for $1 \le i \le 2N+1$.
Let $r_i=\left( r_i^{(1)},\dots,r_i^{(2N+1)} \right)$.
From $r_i A=\lambda_i r_i$ we have
$$r_i^{(j)} + a_j r_i^{(2N+1)} = \lambda_i r_i^{(j+1)} \quad \text{for } 0 \le j \le 2N.$$
Here we have assumed $r_i^{(0)}=0$ for simplicity.
From the last equation we see that $r_i^{(2N+1)} \ne 0$; otherwise the eigenvector $r_i=0$. Thus we may as well assume $r_i^{(2N+1)}=1$.  Recall that $a_{2N+1}=-1$. Then we can easily obtain
\begin{displaymath}
    r_i^{(j)} = - \sum_{k=j}^{2N+1}a_k \lambda_i^{k-j}
\end{displaymath}
for $0 \le j \le 2N$. Therefore, we have
\begin{equation} \label{eq:symmL}
    L=
    \begin{bmatrix}
        \lambda_1^{2N} & \lambda_1^{2N-1} & \lambda_1^{2N-2} & \cdots & 1 \\
        \lambda_2^{2N} & \lambda_2^{2N-1} & \lambda_2^{2N-2} & \cdots & 1 \\
        \lambda_3^{2N} & \lambda_3^{2N-1} & \lambda_3^{2N-2} & \cdots & 1 \\
        \vdots & \vdots & \vdots & & \vdots \\
        \lambda_{2N+1}^{2N} & \lambda_{2N+1}^{2N-1} & \lambda_{2N+1}^{2N-2} & \cdots & 1
    \end{bmatrix}
    \begin{bmatrix}
        1 &&&& \\
        -a_{2N} & 1 &&& \\
        -a_{2N-1} & -a_{2N} & 1 && \\
        \vdots & \vdots & \vdots & \ddots & \\
        -a_1 & -a_2 & -a_3 & \cdots & 1
    \end{bmatrix}.
\end{equation}

With this $L$, we can state our main result of this subsection.

\begin{theorem} \label{thm:eqmombgkstabiii}
    For the EQMOM moment system \cref{eq:EQMOMboltz}, the inequality in the structural stability condition (iii) holds with $A_0 = L^T L$ and $P$ defined in \cref{eq:sourcediag}.
\end{theorem}

\begin{proof}
    According to Theorem 2.1 in \cite{Yong1999}, it suffices to show that at equilibrium states $M$,
    $$K(M) := P^{-T} A_0 P^{-1} = (LP^{-1})^{T} (LP^{-1})$$
    is of the block-diagonal form $\diag(K_1,K_2)$, in which $K_1$ and $K_2$ are $3 \times 3$ and $(2N-2) \times (2N-2)$ matrices, respectively.
    Namely, the first three columns of $LP^{-1}$ are orthogonal to its other columns.
    In what follows all the states $M$ are in equilibrium.

    To show the orthogonality, we compute the $(2N+1)$-matrix $LP^{-1}: = (b_{il})$. From \cref{eq:sourcediag} we see that
    \begin{equation} \label{eq:sourcediaginv}
        P^{-1} =
        \begin{bmatrix}
            I_3 & \\
            \hat{S}_M & I_{2N-2}
        \end{bmatrix}.
    \end{equation}
    This, together with \cref{eq:symmL}, gives
    \begin{displaymath}
        b_{il} = \left \{
        \begin{aligned}
            - \sum_{j=0}^{2N} \chi_j^{l-1} \sum_{k=j+1}^{2N+1} a_k \lambda_i^{k-j-1} \quad &\text{for } 1 \le l \le 3, \\
            -\sum_{k=l}^{2N+1} a_k \lambda_i^{k-l} \quad &\text{for } 4 \le l \le 2N+1.
        \end{aligned}
        \right.
    \end{displaymath}

    The expression above indicates that the last $(2N-2)$ columns of $LP^{-1}$ are linear combinations of $(\lambda_1^{\beta},\dots,\lambda_{2N+1}^{\beta})^{T} \in \mathbb{R}^{2N+1}$ for $0 \le \beta \le 2N-3$.
    Thus it reduces to show that
    $$\sum_{i=1}^{2N+1} b_{il} \lambda_i^{\beta} = 0 \quad \text{for $l=1,2,3$ and $0 \le \beta \le 2N-3$} .$$
    Set $p_k = \sum_{i=1}^{2N+1} \lambda_i^k$. By using the above expression of $b_{il}$ for $1\le l\le 3$, the last equation is equivalent to
    \begin{equation} \label{eq:LPinvinprod}
        \sum_{j=0}^{2N} \chi_j^{l} \sum_{k=j+1}^{2N+1} a_k p_{k-j-1+\beta} =   \sum_{j=-\beta}^{2N-\beta} \chi_{j+\beta}^l \sum_{k=\beta}^{2N-j} a_{j+k+1}p_{k} = 0
    \end{equation}
    for $l=0,1,2$ and $0 \le \beta \le 2N-3$.

    To prove \cref{eq:LPinvinprod}, it suffices to show that
    \begin{equation} \label{eq:LPinvinpdeq}
        \sum_{j=-\beta}^{2N-\beta} \mathcal{H}_{j+\beta} \sum_{k=\beta}^{2N-j} a_{j+k+1} p_{k} = 0, \quad \text{for } 0 \le \beta \le 2N-3,
    \end{equation}
    where $\mathcal{H}_j$ can be replaced by any of $\Delta_j=\Delta_j(U,\theta)$, $\partial_U \Delta_j$ and $\partial_{\theta} \Delta_j$. Indeed, \cref{eq:LPinvinprod} follows immediately from \cref{eq:LPinvinpdeq} and \cref{eq:chidef} which says
    \begin{displaymath}
        \chi_j^l = \frac{\partial}{\partial M_l} (\rho \Delta_j(U,\theta)) = \left( \frac{\partial \rho}{\partial M_l} \right) \Delta_j + \left( \rho \frac{\partial U}{\partial M_l} \right) \partial_U \Delta_{j} + \left( \rho \frac{\partial \theta}{\partial M_l} \right) \partial_{\theta} \Delta_j
    \end{displaymath}
    for  $0 \le j \le 2N$ and $l=0,1,2$.

    Before proceeding, two tools are needed. The first one is Newton's power sum formulas for $p_k$ \cite{EID1968}:
    \begin{subequations} \label{eq:newtonsum}
    \begin{align}
        \sum_{k=0}^{2N+1} a_k p_{k-j-1}  = \sum_{k=-1-j}^{2N-j} a_{j+k+1} p_k &= 0 \quad \text{for} \ j \le -2, \label{eq:newtonsum2} \\
        (2N-j) a_{j+1} + \sum_{k=1}^{2N-j} a_{j+k+1} p_k &= 0 \quad \text{for} \ -1 \le j \le 2N. \label{eq:newtonsum1}
    \end{align}
    \end{subequations}
    The second tool is the following relation
    \begin{equation} \label{eq:ceqeq}
        \left. \mathcal{D}_{\theta} \left( u^k c^{(j)} \right) \right|_{u=U} = 0
    \end{equation}
    for $0 \le k \le 2N$ and $0 \le j \le 2N-k$, where $c^{(j)}$ denotes the $j$th derivative of the characteristic polynomial $c=c(u;W)$ with respect to $u$.
    This relation can be proved as below. \cref{prop:groots} tells $g(u) = (u-U)^{2N+1}$ in equilibrium and therefore $g^{(j)}(U) = 0$ for $0 \le j \le 2N$.
    From \cref{prop:basicOp}(c) we see that $g^{(j)}(u)=\mathcal{D}_{\theta} c^{(j)}$ and thereby $\left. \mathcal{D}_{\theta} c^{(j)} \right|_{u=U} = 0$ for $0 \le j \le 2N$. This is just the case for $k=0$ in \cref{eq:ceqeq}.
    Then using \cref{prop:basicOp}(d) we have
    \begin{displaymath}
        \left. \mathcal{D}_{\theta}  \left(u c^{(j)}\right) \right|_{u=U}
        = U \left. \mathcal{D}_{\theta} c^{(j)} \right|_{u=U} + \theta \left. \mathcal{D}_{\theta} c^{(j+1)} \right|_{u=U} = 0
    \end{displaymath}
    for $j=0,\dots,2N-1$, which validates the case for $k=1$.
    This procedure can be repeated to show  \cref{eq:ceqeq} for other $k \le 2N$.

    With these preparations, we only need to prove \cref{eq:LPinvinpdeq} for the following two cases.

    \textbf{Case I}: $\beta=0$. Noting that $p_0 = 2N+1$, we deduce from \cref{eq:newtonsum1} that $$\sum_{k=0}^{2N-j} a_{k+j+1} p_{k} = (j+1)a_{j+1}.$$
    Thus \cref{eq:LPinvinpdeq} in this case is equivalent to
    \begin{displaymath}
        \sum_{j=0}^{2N} (j+1) a_{j+1} \mathcal{H}_j = 0.
    \end{displaymath}
    When taking $\mathcal{H}_j$ to be $\Delta_j$, $\partial_U \Delta_j$ or $\partial_{\theta} \Delta_j$,
    the left-hand side of the last equation is equivalent to $\left. \mathcal{D}_{\theta} c' \right|_{u=U}$, $\left. \mathcal{D}_{\theta} c'' \right|_{u=U}$ or $\left. \mathcal{D}_{\theta} c''' \right|_{u=U}$, respectively.
    They are all equal to zero due to \cref{eq:ceqeq} and hence \cref{eq:LPinvinpdeq} with $\beta=0$ is proved.

    \textbf{Case II}: $\beta \ge 1$.
    As in Case I, we first simplify the coefficients $\sum_{k=\beta}^{2N-j} a_{j+k+1} p_{k}$ of $\mathcal{H}_{j+\beta}$ in \cref{eq:LPinvinpdeq} by using Newton's power sum formulas \cref{eq:newtonsum1,eq:newtonsum2}. They can be rewritten as
    \begin{displaymath}
    \begin{aligned}
        \sum_{k=-1-j}^{\beta-1} a_{j+k+1} p_{k} + \sum_{k=\beta}^{2N-j} a_{j+k+1} p_{k} &= 0 \quad \text{for } j \le -2, \\
        \left( (2N-j)a_{j+1} + \sum_{k=1}^{\beta-1}a_{j+k+1} p_{k} \right) + \sum_{k=\beta}^{2N-j} a_{j+k+1} p_{k} &= 0 \quad \text{for } j \ge -1.
    \end{aligned}
    \end{displaymath}
    With these two relations, \cref{eq:LPinvinpdeq} is equivalent to
    \begin{equation} \label{eq:ceqfurred}
    \begin{split}
        0
        &= \sum_{j=-\beta}^{2N-\beta} \mathcal{H}_{j+\beta} \sum_{k=\max \{ 1, -1-j \}}^{\beta-1} a_{j+k+1} p_{k} + \sum_{j=-1}^{2N-\beta} (2N-j) \mathcal{H}_{j+\beta} a_{j+1} \\
        &= \sum_{k=1}^{\beta-1} p_{k} \sum_{j=-1-k}^{2N-\beta} \mathcal{H}_{j+\beta} a_{j+k+1} + \sum_{j=-1}^{2N-\beta} (2N-j) \mathcal{H}_{j+\beta} a_{j+1}
    \end{split}
    \end{equation}
    for $1 \le \beta \le 2N-3$.

    \cref{eq:ceqfurred} can be further simplified by using the following relations
    \begin{displaymath}
        \sum_{j=-1-k}^{2N-k} \mathcal{H}_{j+\beta} a_{j+k+1} = 0 \quad \text{for } 1 \le k \le \beta-1.
    \end{displaymath}
    Indeed, replacing $\mathcal{H}_j$ by $\Delta_j$, $\partial_U \Delta_j$ or $\partial_{\theta} \Delta_j$,
    the sum is just $\left. \mathcal{D}_{\theta} \left( u^{\beta-1-k} c \right) \right|_{u=U}$, $\left. \mathcal{D}_{\theta} \left( u^{\beta-1-k} c' \right) \right|_{u=U}$ or $\left. \mathcal{D}_{\theta} \left( u^{\beta-1-k} c'' \right) \right|_{u=U}$,
    respectively.
    They are all equal to zero due to \cref{eq:ceqeq} and the fact that $0 \le \beta-1-k \le \beta-2 \le 2N-5$.
    With the last relation, the first term in the right-hand side of \cref{eq:ceqfurred} is reduced to
    \begin{displaymath}
    \begin{aligned}
        \sum_{k=1}^{\beta-1} p_{k} & \sum_{j=-1-k}^{2N-\beta} \mathcal{H}_{j+\beta} a_{j+k+1}
        = -\sum_{k=1}^{\beta-1} p_{k} \sum_{j=2N-\beta+1}^{2N-k} \mathcal{H}_{j+\beta} a_{j+k+1} \\
        &= -\sum_{j=2N-\beta+1}^{2N-1} \mathcal{H}_{j+\beta} \sum_{k=1}^{2N-j} a_{j+k+1} p_{k}
        = \sum_{j=2N-\beta+1}^{2N-1} \mathcal{H}_{j+\beta} \left[ (2N-j) a_{j+1} \right].
    \end{aligned}
    \end{displaymath}
    The last step resorts again to \cref{eq:newtonsum1} for $j \ge 2N+1-\beta \ge 4$.
    With this, \cref{eq:ceqfurred} is equivalent to
    \begin{equation} \label{eq:ceqredfin}
        \sum_{j=-1}^{2N-1} (2N-j) \mathcal{H}_{j+\beta} a_{j+1} = 0 \quad \text{for } 1 \le \beta \le 2N-3.
    \end{equation}
    This is our final task.

    In \cref{eq:ceqredfin}, we take $\mathcal{H}_j$ to be $\Delta_j$, $\partial_U \Delta_j$ or $\partial_{\theta} \Delta_j$ and arrive at
    \begin{displaymath}
        \left. \mathcal{D}_{\theta} \left( u^{\beta}c' - (2N+1) u^{\beta-1} c \right)^{(k)} \right|_{u=U} = 0
    \end{displaymath}
    for $1 \le \beta \le 2N-3$.
    Here $k=0,1,2$ correspond to $\mathcal{H}_j=\Delta_j, \ \partial_U \Delta_j \text{ or } \partial_{\theta} \Delta_j$, respectively.
    The last relations hold due to \cref{eq:ceqeq} and the linearity of the operator $\mathcal{D}_{\theta}$.
    Hence the orthogonality is validated and the proof is completed.
\end{proof}

\begin{remark}
    It is worth pointing out that the structural stability condition still holds if the collision frequency $\nu=\nu(M)$ in the BGK model depends on $M$, because in equilibrium $S=S(M)=0$ and thus $\partial_M (\nu S) = \nu S_M(M) + S \partial_M \nu = \nu S_M(M)$. Hence all the analyses above are valid.
\end{remark}

\subsection{Shakhov model}
\label{subsec:coupASseqmom}

This subsection is devoted to the EQMOM moment system of the 1-D Boltzmann equation with Shakhov source term \cref{eq:Shakhov}. We first introduce the notation
\begin{displaymath}
    \Delta_j^S =  \Delta_j^S (U,\theta,q) = \frac{1}{\rho} \int_{\mathbb{R}} \xi^j f_S d \xi
\end{displaymath}
with the equilibrium distribution $f_S = f_S(t,x,\xi)$ defined in \cref{eq:fS}. In this situation the moment system has the form \cref{eq:EQMOMboltz} but the source term is different:
\begin{equation*}
    S^{Sh} = S^{Sh}(M) =  \rho (\Delta_0^S, \Delta_1^S, \dots, \Delta_{2N}^S)^{T} - M.
\end{equation*}

We need to investigate whether this source term satisfies the structural stability condition (i) \& (iii). For this purpose, some basic properties of $\Delta_j^S$ are required. A direct calculation shows that $\Delta_0^S = \Delta_0 = 1$, $\Delta_1^S = \Delta_1 = U$ and $\Delta_2^S = \Delta_2 = U^2+\theta$. Moreover, for $j \ge 3$ we have
\begin{equation} \label{eq:deltaS}
\begin{split}
    \rho \Delta_j^S - \rho \Delta_j
    &= \frac{q(1-Pr)}{3 \theta^2} \int_{\mathbb{R}} \xi^j (\xi-U) \left( \frac{(\xi-U)^2}{\theta}-3 \right) f_{eq} d\xi \\
    &= \binom{j}{3} (1-Pr) (2q) \Delta_{j-3}
      = \binom{j}{3} (1-Pr) (M_3 - \rho \Delta_3) \Delta_{j-3}.
\end{split}
\end{equation}
Here \cref{prop:basicDel}(a) is used for the integration and the last step is due to the definition of $q$.

As in \Cref{subsec:eqmomequistate}, the equilibrium state $W$ needs to be determined. From $S^{Sh}(M)=0$ we see that $\rho \Delta_j^S = M_j$ for $0 \le j \le 2N$.
With \cref{eq:deltaS}, the equation $\rho \Delta_3^S = M_3$ clearly implies that $M_3=\rho \Delta_3$ for $Pr \ne 1$.
Therefore, we have $\rho \Delta_j^S = \rho \Delta_j$ for any $0 \le j \le 2N$ and the equilibrium manifold $\mathcal{E}$ is determined by $M_j = \rho \Delta_j$ for any $0\le j \le 2N$. This is exactly the same as that of the BGK model, which has already been determined in \cref{thm:equisolu} to be $W \in \Omega_W^{eq}$.

At equilibrium, the Jacobian matrix of $S^{Sh}$ can be computed with \cref{eq:deltaS}:
\begin{equation} \label{eq:sourcejacs}
    S_M^{Sh}: = \left. \frac{\partial S^{Sh}}{\partial M} \right|_{S^{Sh}(M)=0} =
    \left( I_{2N+1} - (1-Pr) \sum_{i=3}^{2N} \binom{i}{3} \Delta_{i-3} E_{(i+1),4} \right) S_M,
\end{equation}
where $S_M$ is the Jacobian matrix \cref{eq:sourcejac} for the BGK model and the $(2N+1)$-matrix $E_{ij} = \left( e_{ij} \right)$ with $e_{ij}=1$ and all the other entities being zero.
$S_M^{Sh}$ is diagonalizable by an invertible matrix $P^S$ such that $P^S S_M^{Sh} = -\diag(0,0,0,Pr,1,\dots,1) P^S$, and
\begin{equation} \label{eq:sourcediags}
    \left(P^S \right)^{-1} =
    P^{-1} + \sum_{i=4}^{2N} \binom{i}{3} \Delta_{i-3} E_{(i+1),4},
\end{equation}
where $P^{-1}$ is defined in \cref{eq:sourcediaginv}.
Hence the structural stability condition (i) is justified.

For Condition (iii), we take the same symmetrizer $A_0=L^{T} L$ as that for the BGK model. This is reasonable since the equilibrium state is the same. It then suffices to show that the first three columns of $L \left( P^S \right)^{-1}$ are orthogonal to its other columns in equilbrium.
From \cref{eq:sourcediags} we see that the only difference between $L \left( P^S \right)^{-1}$ and $LP^{-1}$ is the fourth column.
For $L \left( P^S \right)^{-1}$, its fourth column is a linear combination of the last $(2N-2)$ columns of $L$, while the last $(2N-2)$ columns of $LP^{-1}$ are exactly those of $L$.
Since the first three columns of $LP^{-1}$ are orthogonal to its other columns, the fourth column of $L \left( P^S \right)^{-1}$ is also orthogonal to its first three columns. This has validated Condition (iii).
In this way, we have the main result of this subsection:

\begin{theorem} \label{thm:eqmomshakstab}
    For the 1-D Boltzmann equation with the Shakhov model, the EQMOM moment system satisfies the structural stability condition.
\end{theorem}

\section{Conclusions}
\label{sec:conclusions}

This paper presents a rigorous stability analysis of the quadrature based moment methods (QBMM) for the Boltzmann equation.
To figure out a road map for more general cases, only the spatial one-dimensional (1-D) Boltzmann equation with hypothetical collisions (BGK or Shakhov type) is considered here.
In the QBMM, the distribution function $f$ is approximated with a linear combination of $N$ ($N \ge 1$) $\delta$-functions with unknown centers or their Gaussian approximations with unknown variance and centers (named QMOM or EQMOM, respectively).
For QMOM, we show purely analytically that the resulting moment systems of first-order PDEs are not strongly hyperbolic for any $N$.
Furthermore, we prove that the moment systems produced by the Gaussian EQMOM are strictly hyperbolic, when the variance is positive, and preserve the dissipation property of the kinetic equation. As a step in the proof, we also determine the equilibrium manifold that lies on the boundary of the state space for the parameters $(w_i, u_i, \sigma^2)$ ($1\le i \le N$). These conclusions explain why the EQMOM gives reasonable numerical results while QMOM does not.

The proofs are quite technical and involve detailed analyses of the characteristic polynomial of the coefficient matrices. They offer a guideline to investigate the multidimensional cases with multiple nodes, which is underway.

\appendix
\section{Injectivity of EQMOM}
\label{sec:injection}
In this appendix we show
\begin{proposition}
    For EQMOM, the map $M = \mathcal{M} (W)$ in \cref{eq:EQMOMmo} is injective for $W \in \Omega_W^{open}$ defined in \cref{eq:domainWopen}.
\end{proposition}

\begin{proof}
    It suffices to demonstrate that the Jacobian matrix $\frac{\partial \mathcal{M}}{\partial W}$ in \cref{eq:jacoeqmom} is invertible for $W \in \Omega_W^{open}$. In fact, we can show that
    \begin{displaymath}
        \det{\left( \frac{\partial \mathcal{M}}{\partial W} \right)} = \left( \prod_{i=1}^N w_i \right) \cdot \left( \sum_{i=1}^N w_i \prod_{\substack{j=1 \\ j\ne i}}^N (u_i-u_j)^2 \right) \cdot \prod_{1 \le i < j \le N} (u_i-u_j)^4
    \end{displaymath}
    for multiple nodes $N \ge 2$.
    To this end, we set $\mathcal{F}(u) = \left( \Delta_0(u,\sigma^2), \Delta_1(u,\sigma^2), \cdots, \Delta_{2N}(u,\sigma^2) \right)^T$ and see from \cref{eq:jacoeqmom} that
    \begin{equation*}
        \begin{aligned}
            \det \left(\frac{\partial \mathcal{M}}{\partial W} \right)
            &= \left( \prod_{i=1}^N w_i \right) \cdot \det \left( \mathcal{F}(u_1), \mathcal{F}'(u_1), \cdots, \mathcal{F}(u_N), \mathcal{F}'(u_N), \frac{1}{2} \sum_{i=1}^N w_i \mathcal{F}''(u_i) \right) \\
            &= \left( \prod_{i=1}^N w_i \right) \cdot \sum_{i=1}^N w_i \det \left( \mathcal{F}(u_1), \mathcal{F}'(u_1), \cdots, \mathcal{F}(u_N), \mathcal{F}'(u_N), \frac{1}{2} \mathcal{F}''(u_i) \right).
        \end{aligned}
    \end{equation*}

    Denote
    \begin{equation*}
        f_N (u_1,u_2,\dots,u_N;\sigma) = \det \left( \mathcal{F}(u_1), \mathcal{F}'(u_1), \cdots, \mathcal{F}(u_N), \mathcal{F}'(u_N), \frac{1}{2} \mathcal{F}''(u_1) \right).
    \end{equation*}
    And we see that for each $1\le i \le N$,
    \begin{equation*}
        \det \left( \mathcal{F}(u_1), \mathcal{F}'(u_1), \cdots, \mathcal{F}(u_N), \mathcal{F}'(u_N), \frac{1}{2} \mathcal{F}''(u_i) \right) = f_N (u_i,u_2,\dots,u_{i-1}, u_1,u_{i+1}, \dots,u_N;\sigma),
    \end{equation*}
    and thereby
    \begin{equation*}
        \det \left(\frac{\partial \mathcal{M}}{\partial W} \right)
        = \left( \prod_{i=1}^N w_i \right) \cdot \sum_{i=1}^N w_i
            f_N (u_i, u_2, \dots, u_{i-1}, u_1, u_{i+1}, \dots, u_N;\sigma).
    \end{equation*}
    Thus, it remains to show
    \begin{equation} \label{eq:detfeqmom}
        f_N(u_1,u_2,\dots,u_N;\sigma) = C(N) \prod_{j=2}^N (u_j-u_1)^2 \cdot \prod_{1 \le i < j \le N} (u_i-u_j)^4
    \end{equation}
    and
    \begin{equation} \label{eq:detfcoef}
        C(N) = 1 \quad \text{for } N \ge 2.
    \end{equation}

    Note that $f_N (u_1,\dots,u_N;\sigma)$ is a homogeneous polynomial of $u_1,\dots,u_N, \sigma$ with degree $2(N^2-1)$. This can be seen from the definition of determinant and the fact that $\Delta_j (u,\sigma^2)$ is a homogeneous polynomial of $u$ and $\sigma$ with degree $j$ (see \cref{prop:basicDel}(a)).
    On the other hand, the right-hand side of \cref{eq:detfeqmom} is also a homogeneous polynomial of $u_1,\dots,u_N, \sigma$ with degree $2(N^2-1)$.
    Thus, to prove \cref{eq:detfeqmom}, we need to show that $(u_j-u_1)^6$ and $(u_j-u_i)^4$ are factors of $f(u_1,\dots, u_N;\sigma)$ for any $2 \le j \ne i \le N$.

    From the definition of $f_N=f_N(u_1,\dots,u_N;\sigma)$, it is not difficult to compute that for $j\ne 1$,
    \begin{equation*}
    \begin{aligned}
        \partial_{u_j} f_N
        =& \det \left( \mathcal{F}(u_1), \mathcal{F}'(u_1), \cdots, \mathcal{F}(u_j), \mathcal{F}''(u_j), \cdots, \frac{1}{2} \mathcal{F}''(u_1) \right), \\
        \partial_{u_j}^2 f_N
        =& \det \left(\cdots, \mathcal{F}'(u_j), \mathcal{F}''(u_j), \cdots \right) + \det \left( \cdots, \mathcal{F}(u_j), \mathcal{F}'''(u_j), \cdots \right), \\
        \partial_{u_j}^3 f_N
        =& \det \left( \cdots, \mathcal{F}(u_j), \mathcal{F}^{(4)}(u_j), \cdots \right) + 2 \det \left( \cdots, \mathcal{F}'(u_j), \mathcal{F}'''(u_j), \cdots \right), \\
        \partial_{u_j}^4 f_N
        =& \det \left( \cdots, \mathcal{F}(u_j), \mathcal{F}^{(5)}(u_j), \cdots \right) + 3 \det \left( \cdots, \mathcal{F}'(u_j), \mathcal{F}^{(4)}(u_j), \cdots \right) \\
        &+ 2 \det \left( \cdots, \mathcal{F}''(u_j), \mathcal{F}'''(u_j), \cdots \right), \\
        \partial_{u_j}^5 f_N
        =& \det \left( \cdots, \mathcal{F}(u_j), \mathcal{F}^{(6)}(u_j), \cdots \right) + 4 \det \left( \cdots, \mathcal{F}'(u_j), \mathcal{F}^{(5)}(u_j), \cdots \right) \\
        &+ 5 \det \left( \cdots, \mathcal{F}''(u_j), \mathcal{F}^{(4)}(u_j), \cdots \right).
    \end{aligned}
    \end{equation*}
    Thus it follows that for $j \ne 1$,
    \begin{equation*}
        \left. f_N \right |_{u_j = u_i} =0, \quad
        \left. \partial_{u_j} f_N \right |_{u_j = u_i} =0, \quad
        \left. \partial_{u_j}^2 f_N \right |_{u_j = u_i} =0, \quad
        \left. \partial_{u_j}^3 f_N \right |_{u_j = u_i} =0
    \end{equation*}
    for any $1\le i \ne j \le N$ and that
    \begin{equation*}
        \left. \partial_{u_j}^4 f_N \right |_{u_j = u_1} =0, \quad
        \left. \partial_{u_j}^5 f_N \right |_{u_j = u_1} =0.
    \end{equation*}
    This justifies \cref{eq:detfeqmom} and $C(N)$ is a constant.

    Then all we need is to prove \cref{eq:detfcoef}. A direct calculation for $N=2$ indicates that $f_2 (u_1,u_2;\sigma) = (u_1-u_2)^4$ and thus $C(2)=1$. For $N>2$, we deduce from \cref{eq:detfeqmom} that the leading coefficient of $u_N$ (with degree $(4N-2)$) is
    \begin{equation*}
        C(N) \prod_{j=2}^{N-1}(u_j-u_1)^2 \cdot \prod_{1 \le i < j \le N-1} (u_i-u_j)^4
    \end{equation*}
    On the other hand, the determinant definition of $f_N(u_1,\dots,u_N;\sigma)$ implies that the leading term of $u_N$ (with degree $(4N-2)$) is included in the following part:
    \begin{equation*}
        f_{N-1}(u_1,\dots,u_{N-1};\sigma) \times
        \det
        \begin{bmatrix}
            \Delta_{2N-1}(u_N) & (2N-1)\Delta_{2N-2}(u_N) \\
            \Delta_{2N}  (u_N) & 2N\Delta_{2N-1}(u_N)
        \end{bmatrix}.
    \end{equation*}
    Thus, using \cref{eq:detfeqmom}, the leading coefficient of $u_N$ is
    \begin{equation*}
        f_{N-1}(u_1,\dots,u_{N-1};\sigma) = C(N-1) \prod_{j=2}^{N-1}(u_j-u_1)^2 \cdot \prod_{1 \le i < j \le N-1} (u_i-u_j)^4.
    \end{equation*}
    By equating the leading coefficients of $u_N$ in the above two expressions, we see immediately that $C(N)=C(N-1)$ for $N > 2$. Since $C(2)=1$, this justifies \cref{eq:detfcoef} and hence completes the proof.
\end{proof}

\bibliographystyle{siamplain}
\bibliography{references}
\end{document}


\maketitle

\section{1-D Riemann problem}

In this part we report our numerical experiments by using the two-node QMOM and EQMOM to solve the 1-D Riemann problem with the BGK model.
Assume that the collision frequency is proportional to the density, i.e. $\nu= \kappa \rho$. Three values of $\kappa$, $0,10,\infty$ are considered.
When $\kappa=0$, $f$ is a traveling wave and the analytical solution is given in \cite{Chalons2017}. When $\kappa=\infty$, the analytical solution of the Euler equation can be found in \cite{Toro2009}. In contrast, the direct simulation Monte Carlo (DSMC) result for $\kappa=10$ is used for a comparison, in which the BGK collision term is simulated according to \cite{Macro2002}.

In the simulation, the initial values are \cite{Chalons2017}: $\rho(0,x)=1$ and $\theta(0,x) = \frac{1}{3}$ for any $x \in \mathbb{R}$; $U(0,x<0) = 1$ and $U(0,x>0)=-1$. The moment inversion algorithms and spatial fluxes are treated as in \cite{Chalons2017,Fox2008}. The 1-D computational domain of $-1<x<1$ is discretized into 1000 uniform cells.
The time step is chosen so that the CFL number is less than 0.5.

\begin{figure}[htbp]
  \centering
  \label{fig:kpa0}\includegraphics[width=5in]{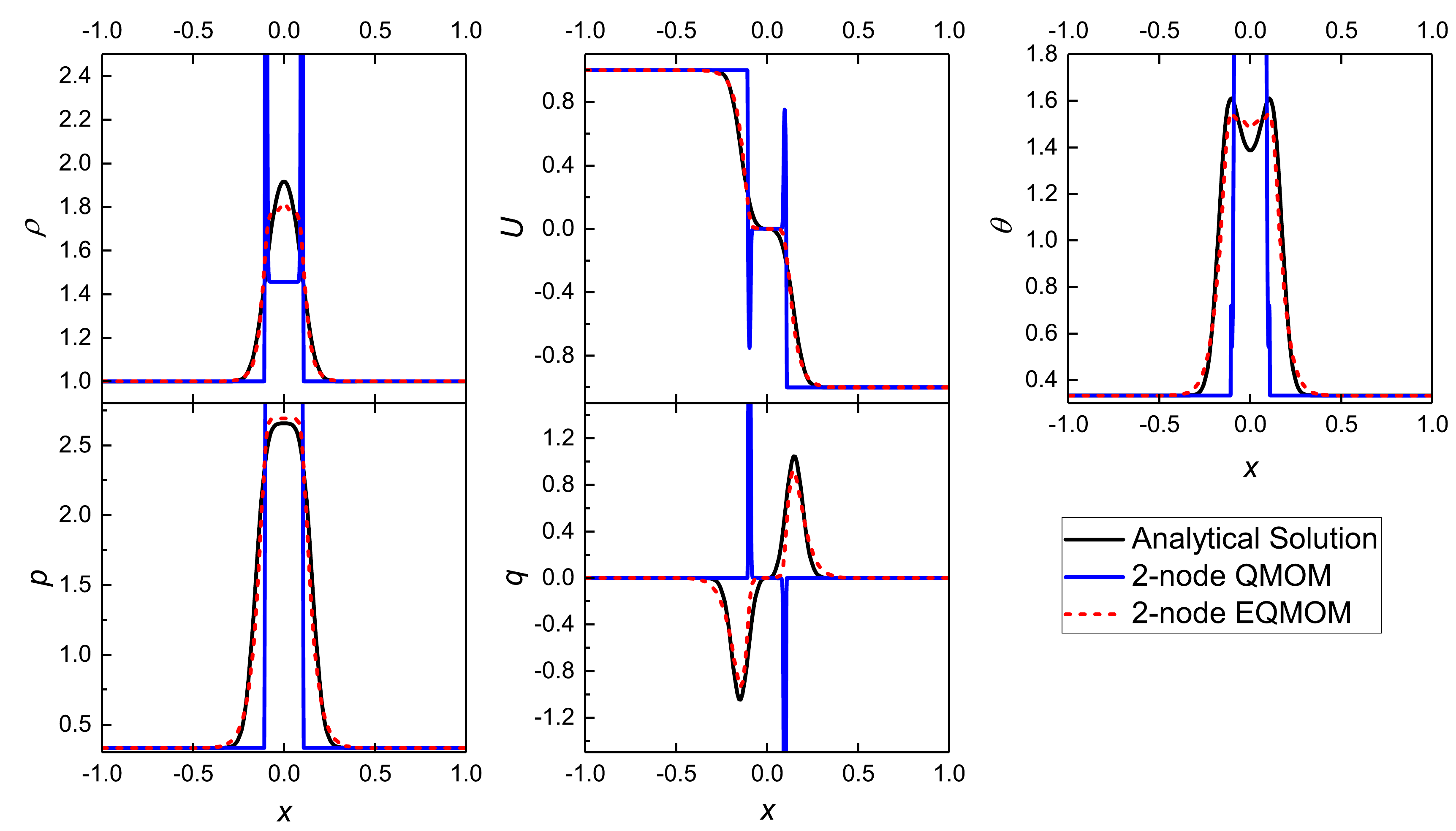}
  \caption{The macroscopic parameters at $t=0.1$ for $\kappa=0$.}
  \label{fig:numfigkpa0}
\end{figure}

\begin{figure}[htbp]
  \centering
  \label{fig:kpa10}\includegraphics[width=5in]{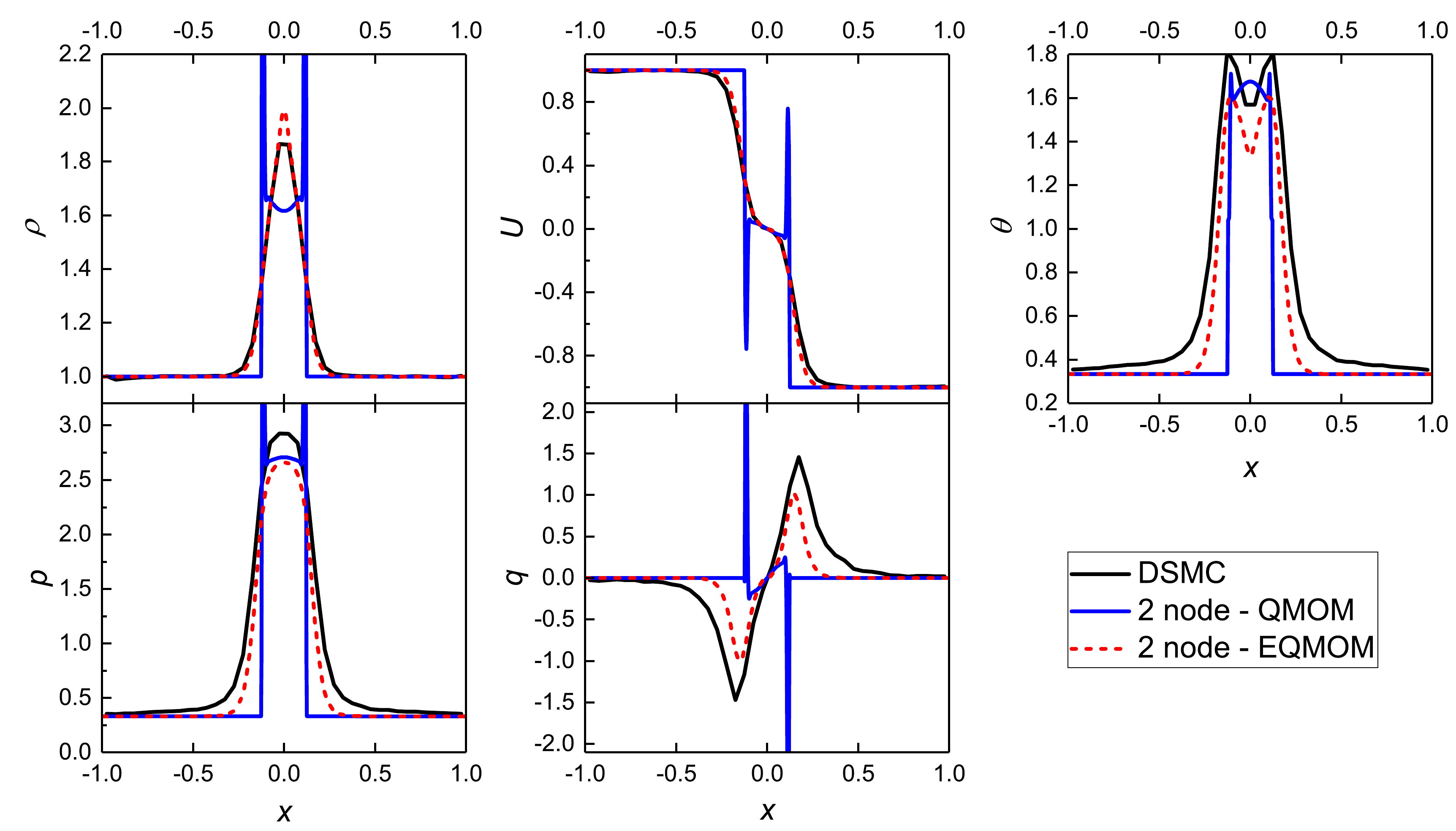}
  \caption{The macroscopic parameters at $t=0.1$ for $\kappa=10$.}
  \label{fig:numfigkpa10}
\end{figure}

\begin{figure}[htbp]
  \centering
  \label{fig:kpainf}\includegraphics[width=5in]{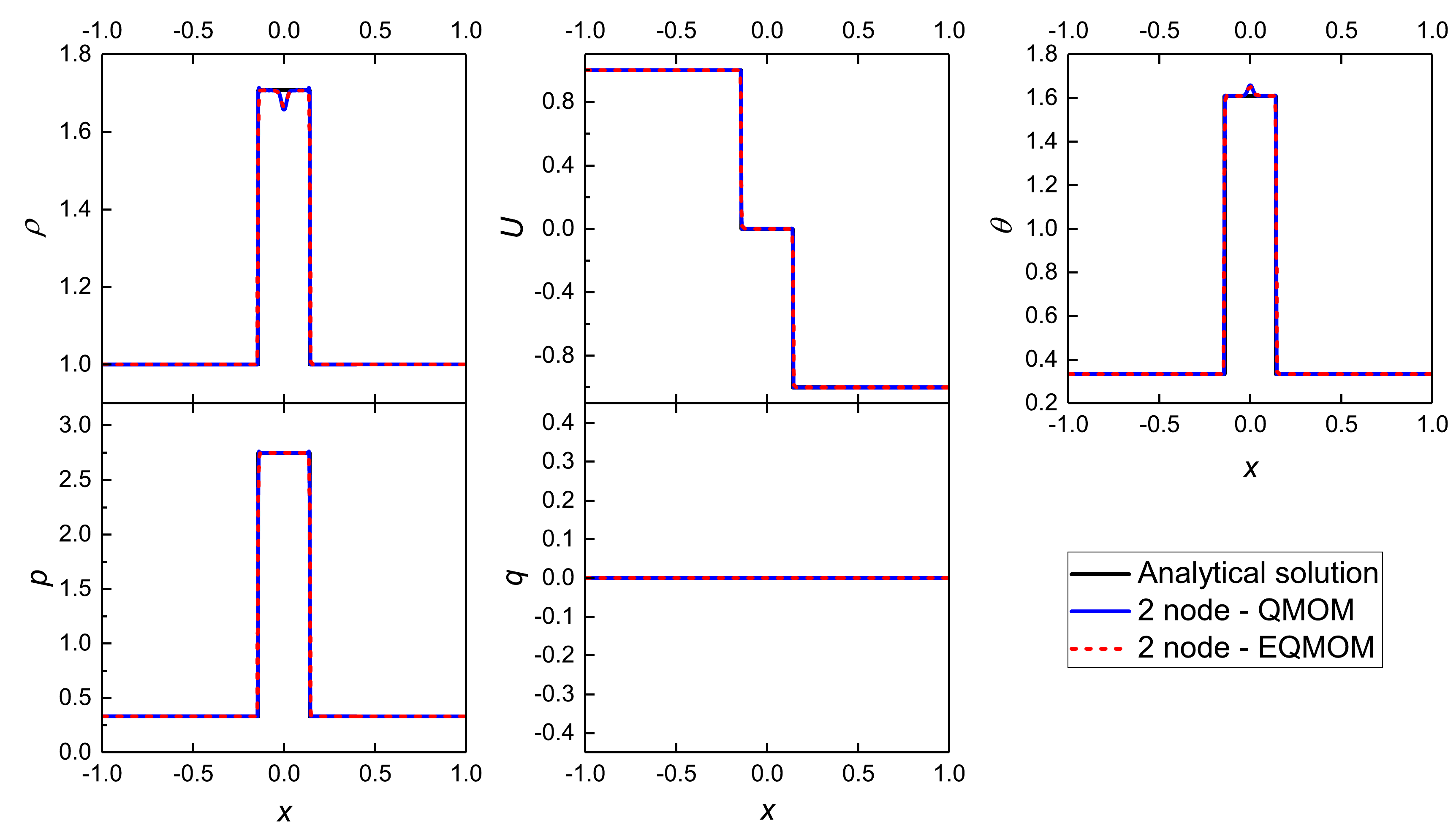}
  \caption{The macroscopic parameters at $t=0.1$ for $\kappa=\infty$.}
  \label{fig:numfigkpainf}
\end{figure}

The macroscopic parameters defined in \cref{eq:mo2macro} at $t=0.1$ are shown in \cref{fig:numfigkpa0,fig:numfigkpa10,fig:numfigkpainf} for the three values of $\kappa$. When the collision rate is zero or finite, QMOM results in a huge unphysical $\delta$-shock, whereas the EQMOM result is close to either the analytical solution or the DSMC results.
Such numerical results demonstrate again that EQMOM gives reasonable and much better results. This motivates us to study the two methods analytically.

\bibliographystyle{siamplain}
\bibliography{references}